\documentclass{elsarticle}

\usepackage{amsfonts}
\usepackage{amsbsy}
\usepackage{graphicx}
\usepackage{amsmath}
\usepackage{xcolor}
\usepackage{ulem}

\newcommand{\PF}[1]{\textcolor{black}{#1}}

\newdefinition{rmk}{Remark}

\title{High resolution compact implicit numerical scheme for conservation laws\tnoteref{vega}}
\tnotetext[vega]{This work was supported by grants VEGA 1/0709/19 and APVV-19-0460}

\author[1]{Peter Frolkovi\v{c}\corref{cor}}
\ead{peter.frolkovic@stuba.sk}

\author[1]{Michal \v{Z}erav\'y}
\ead{michal.zeravy@stuba.sk}

\cortext[cor]{Corresponding author}
\address[1]{Department of Mathematics and Descriptive Geometry, Faculty of Civil Engineering, STU Bratislava, Slovakia}

\begin{document}

% REQUIRED
\begin{abstract}
We present a novel \PF{implicit} scheme for the numerical solution of time-dependent conservation laws. The core idea of the presented method \PF{is to exploit and approximate the mixed spatial-temporal} derivative of the solution that occurs naturally when deriving \PF{some second order accurate schemes in time}. Such an approach is introduced in the context of the Lax-Wendroff (or Cauchy-Kowalevski) procedure when the \PF{second} time derivative is not completely replaced by space derivatives using the PDE, but the mixed derivative is kept. If approximated in a suitable way, \PF{the resulting compact implicit scheme produces} algebraic systems that have a more convenient structure than the systems derived by fully implicit schemes. We derive a high resolution TVD form of the implicit scheme for some representative hyperbolic equations in the one-dimensional case, including illustrative numerical experiments.
\end{abstract}

% REQUIRED

\maketitle
% REQUIRED
% \begin{AMS}
%   65M06, 65M22, 
% \end{AMS}
%\newpageafter{abstract}

\section{Introduction}

One of the most straightforward ways to solve numerically the time dependent hyperbolic partial differential equations (PDEs) is to use a Method Of Lines (MOL) when the spatial and the temporal parts of the PDE are discretized separately. Typically, a discretization of the space is performed in the first step that approximates the PDE by a system of ordinary differential equations (ODEs), for which, afterwards, a chosen numerical integration is used. 

Quite naturally, the first candidates to obtain numerical solutions of ODEs are explicit methods that deliver numerical solutions without the need to solve any algebraic system of equations. Such approximations require a careful choice of discretization steps not only due to accuracy requirements but also for stability reasons. The second requirement is specific to numerical methods, and, if not followed, an unstable behavior of the numerical solution can occur even for well-posed problems. For many types of PDEs and the methods for their numerical solutions, such stability requirements are well understood, typically formulated in the form of a stability condition on the choice of time steps. If the known stability restriction does not limit the choice of discretization steps more than the accuracy requirement, the usage of explicit methods is well justified.

Nevertheless, in several cases such restrictions can be too demanding or simply hard to follow, and, consequently, implicit time discretization methods are also considered to avoid the strict stability restrictions. Such methods are receiving increasing attention, especially for the models that describe several dynamic processes with different characteristic speeds, of which only those with slow or moderate speed are of practical interest. In such cases, at least the terms in PDEs related to the processes with the fastest speed are treated implicitly. As a prominent example (that we do not treat here), one can mention the so-called ``all Mach number'' solvers of Euler equations that are treated typically with semi-implicit methods \cite{park_multiple_2005,cordier_asymptotic-preserving_2012,boscarino_high_2016,tavelli_pressure-based_2017} or Implicit-Explicit (IMEX) methods \cite{pareschi2005implicit,boscarino2009class}, see for some recent developments the list of publications (with no intention to be comprehensive)  \cite{abbate_asymptotic-preserving_2019,avgerinos_linearly_2019,boscheri_second_2020,zeifang2020novel,boscheri_high_2021,busto_semi-implicit_2021,boscarino_high_2022} and the references therein. 

The implicit methods do not explicitly define the values of numerical solutions, but instead they provide algebraic systems to be fulfilled by the discrete numerical values. Therefore, to find the numerical solution, a system of algebraic equations must be solved. Clearly, this is the main price to be paid by the usage of implicit methods that must be well justified, especially if the algebraic systems are nonlinear. Consequently, implicit methods that simplify the task of (nonlinear) algebraic solvers are in demand, and that is the main motivation of our study.

To fulfill such a demand, we are interested in time discretization methods that offer a coupled treatment of the temporal and the spatial discretization methods. In this framework, we mention the explicit methods based on Taylor series expansions, such as the Lax-Wendroff procedure 
\cite{qiu_finite_2003,zorio_approximate_2017,carrillo2019compact,carrillo2021lax}, the local time-space DG or FV discretizations \cite{dumbser_finite_2008,dumbser_unified_2008,gassner_explicit_2011,gaburro_posteriori_2021}, and the multi-derivative schemes \cite{seal_high-order_2014}. Regarding the methods based on the Taylor series, several authors have recognized that it can be advantageous to discretize each term in the Taylor expansion with different spatial discretizations \cite{qiu_finite_2003, seal_high-order_2014, tsai_two-derivative_2014, li_two-stage_2019,zeifang_two-derivative_2021} and, eventually, to use the mixed spatial-temporal derivatives in the Lax-Wendroff procedure to obtain the schemes with a more compact stencil \cite{zorio_approximate_2017,carrillo2019compact, carrillo2021lax}. A natural question arises if such an approach can be used to derive some compact implicit schemes for hyperbolic problems.

Such possibility has been recognized in the context of the nonconservative advection equation used in the level set methods in \cite{frolkovic2016semi, frolkovic2018semi}, where a special coupling of temporal and spatial discretization is used to derive a second order accurate unconditionally stable numerical schemes producing simpler algebraic systems than the ones obtained with fully implicit schemes. This type of scheme was introduced under the abbreviation IIOE (``Inflow-Implicit/Outflow-Explicit'') finite volume method for the scalar advection equation in \cite{mikula2011inflow,mikula2014inflow} and later successfully applied in, e.g., \cite{frolkovic2015semi, hahn2019iterative, ibolya2020numerical}. 
As the notation ``Inflow-Implicit and Outflow-Explicit'' suggests, the spatial and temporal discretizations are aware of each other, and they are used in a coupled way.  

Opposite to level set equations where the solution is supposed to be continuous, the hyperbolic problems allow for discontinuous solutions. For this type of problems, the numerical scheme must be conservative to approximate correctly the movement of shock waves, and the scheme must be nonlinear even for linear problems if higher order scheme with non-oscillatory numerical solutions is required. In this context, we introduce a novel compact implicit method for some representative models of hyperbolic equations. Our method can be viewed as a high resolution extension of the first order implicit method presented in \cite{lozano_implicit_2021} where a fractional step method is conveniently combined with a fast sweeping method. Some similarities of our method can be found when compared to the implicit method described in \cite{duraisamy_implicit_2007} using a non-oscillatory reconstruction in space and time. The idea of Taylor series using the mixed spatial-temporal derivatives for the solution of hyperbolic problems in the framework of explicit methods can be found in \cite{zorio_approximate_2017, carrillo2019compact, carrillo2021lax}, similar methods for the solution of ODEs are described in \cite{baeza_reprint_2018} including implicit methods in \cite{baeza_approximate_2020}. Our scheme can be viewed as a compact scheme, as it reduces the  stencil of numerical solution compared to some standard finite difference scheme similarly to approaches in, e.g., \cite{carrillo2019compact,clain_compact_2022}.  Our attempt in this paper is directed at offering an implicit ``black-box'' solver for some hyperbolic problems similarly to \cite{arbogast2020third, puppo_quinpi_2022,gomez-bueno_implicit_2022} with an aim of introducing the motivation and the basic ideas of the high resolution compact implicit scheme.

The paper is organized as follows. In Section \ref{sec1} we introduce the basic notations for the scalar case including the first order accurate numerical scheme. In Section \ref{sec-odvodenie} we derive the second order accurate compact implicit scheme, and we introduce the high resolution form of numerical fluxes. In Section \ref{sec2} we present and motivate the definitions of parameters in the high resolution scheme with the details for the linear advection equation. In Section \ref{sec3} we present algorithmic details of the high resolution method for the nonlinear scalar case, and in Section \ref{sec4} we give additional details for hyperbolic systems. The Section \ref{sec5} on numerical experiments presents several test examples that illustrate the properties of the method. Finally, in Section \ref{sec-conc} we make some concluding remarks. 

\section{Scalar conservation laws}
\label{sec1}

In this section, we aim to solve numerically the scalar nonlinear hyperbolic equation written in the form
\begin{equation}
    \label{cl}
    u_t + f(u)_x = 0 \,, \quad u(x,0) = u^0(x) \,, \,\, x \in R \,, \, t > 0 \,,
\end{equation}
where $u=u(x,t)$ is the unknown function with initial values prescribed by a given function $u^0$, and $f$ is a given flux function. 

To discretize (\ref{cl}), we follow the approach of conservative finite difference methods as described, e.g., in \cite{shu_essentially_1998,qiu_finite_2003,lozano_implicit_2021}. For this purpose, we use the standard notation for the grid nodes $x_i$, $i=0,1,\ldots,I$ with a uniform step $h \equiv x_i-x_{i-1}$ and the discrete times $0=t^0<t^1<\ldots$ with $\tau \equiv t^{n+1}-t^n$, $n=0,1,\ldots,N$, where the integers $I$ and $N$ are given. \PF{Note that the nodes $x_0$ and $x_I$ are the boundary nodes where we have to specify boundary conditions later.} Our aim is to find the approximations $u_i^{n+1} \approx u(x_i,t^{n+1})$. To do so, we follow the standard form of conservative schemes,
\begin{equation}
%\nonumber
    \label{step0}
    u_i^{n+1} + \frac{\tau}{h}\left( F_{i+1/2} - F_{i-1/2} \right) = u_i^n \,,
\end{equation}
where the numerical fluxes \PF{$F_{i+1/2} \approx f(u(x_{i+1/2},t^{n+1/2}))$ will be specified later. Here, $x_{i+1/2}=x_i+h/2 $ and $t^{n+1/2}=t^n + \tau/2$. }

\PF{To derive our method}, we use the approach of the fractional step method presented in \cite{lozano_implicit_2021}. \PF{Firstly, one chooses a flux splitting where} the flux function $f$ is split into the sum of two functions having nonnegative and nonpositive derivatives,
\begin{equation}
    \label{split}
    f = f^+ + f^- \,, \quad \frac{df^+}{du} \ge 0 \,,\,\, \frac{df^-}{du}\le 0 \,, \,\, u \in R \,.
\end{equation}
\PF{A typical choice to obtain (\ref{split}) is the Lax-Friedrichs flux vector splitting} \cite{shu_essentially_1998},
\begin{equation}
    \label{lf}
    f^+(u) = \frac{1}{2} \left( f(u) + \alpha u\right) \,, \quad f^-(u) = \frac{1}{2} \left( f(u) - \alpha u\right) \,,
\end{equation}
where the parameter $\alpha$ is fixed at the maximum value of $|f'(u)|$ over the \PF{considered} values of $u$. \PF{Note that other definitions of the splitting in (\ref{split}) can be used as we do later for Burgers' equation.}
 
Having the splitting (\ref{split}), \PF{one can split analogously the numerical fluxes written formally as}
\begin{equation}
    \label{splitF}
F_{i+1/2} = F_{i+1/2}^{+} + F_{i+1/2}^{-} \,.
\end{equation}
\PF{ Following \cite{lozano_implicit_2021}, we replace (\ref{step0}) using}
the simplest variant of the fractional step method combined with the fast sweeping method that consists of two partial steps. The first step is given by solving the algebraic equations \PF{in the prescribed order,}
\begin{equation}
\label{step1}
u_i^{n+1/2} + \frac{\tau}{h} F_{i+1/2}^{+} = u_i^{n} + \frac{\tau}{h} F_{i-1/2}^{+} \,, \,\, i=1,2,\ldots,I \,,
\end{equation}
where the numerical \PF{fluxes $F_{1/2}^{+}$ and  $F_{I+1/2}^{+}$ shall be determined from boundary conditions, see later some discussions for particular schemes}.
The second step is given by solving the algebraic equations \PF{in the reverse order}
\begin{equation}
\label{step2}
u_i^{n+1} - \frac{\tau}{h} F_{i-1/2}^{-} = u_i^{n+1/2} - \frac{\tau}{h} F_{i+1/2}^{-} \,, \,\, i=I-1,I-2,\ldots,0 \,,
\end{equation}
where now the numerical fluxes $F_{-1/2}^-$ and $F_{I-1/2}^-$ have to be determined from the boundary conditions. 

\PF{The numerical method (\ref{step1}) - (\ref{step2}) has to be completed by some definitions of} the numerical fluxes $F_{i+1/2}^{+}$ in (\ref{step1}) and $F_{i-1/2}^{-}$ in (\ref{step2}). In \cite{lozano_implicit_2021} they are given by the implicit first order accurate upwind scheme,
\begin{equation}
\label{F}
F_{i+1/2}^{+} = f_{i}^{+,n+1} := f^+(u_{i}^{n+1}) \,, \quad
F_{i-1/2}^{-} = f_{i}^{-,n+1} := f^-(u_{i}^{n+1}) \,.
\end{equation}
\PF{In the case of linear advection equation, the method (\ref{step1}) - (\ref{F}) is unconditionally stable having no restriction on the choice of time step using the von Neumann stability analysis \cite{frolkovic_semi-implicit_2021}. Concerning nonlinear hyperbolic problems, a quite extensive numerical evidence of stable behavior for (\ref{step1}) - (\ref{F}) is provided in \cite{lozano_implicit_2021}.}

\PF{Concerning the boundary conditions, we consider here only the options 
\begin{equation}
\label{bc}
F^+_{1/2}=f^+(u_0^{n+1}) \,, \quad F^-_{I-1/2}=f^-(u_I^{n+1}) \,.
\end{equation}
The values $u_0^{n+1}$ (and analogously for $u_I^{n+1}$) are either given by inflow boundary conditions or extrapolated, e.g., $u_0^{n+1}=u_1^{n+1}$ in the case of outflow boundary condition; see also \cite{lozano_implicit_2021}. Note that the fluxes $F^+_{I+1/2}$ and $F^-_{-1/2}$ are well defined in (\ref{F}).}

The important advantage of the proposed method (\ref{step1}) - (\ref{bc}) is that each algebraic equation contains only one unknown $u_i^{n+1}$. The main disadvantage is the low order accuracy that we aim to improve here. Note that in our numerical experiments we use the fractional step method in the first order accurate form (\ref{step1}) - (\ref{step2}), for some higher order extensions see a discussion in \cite{lozano_implicit_2021}.

\section{Compact implicit second order accurate scheme}
\label{sec-odvodenie}

To put our method in the context of previous works and to show its novelty, we derive it in the context of Taylor methods and the Lax-Wendroff procedure \cite{harten_uniformly_1997,shu_essentially_1998,qiu_finite_2003,leveque_finite_2004,toro_riemann_2009,baeza_reprint_2018,frolkovic2018semi,baeza_approximate_2020}. For simplicity, we consider $f^+ \equiv f$ and $f^- \equiv 0$, so we can derive our numerical scheme directly using the one step form in (\ref{step0}) without two fractional steps (\ref{step1}) - (\ref{step2}). 

For brevity, let $U := u(x_i,t^{n+1})$ and similarly for the derivatives of $u$. Moreover, we suppose that the values at the time level $t^n$ are known and exact, so $u_i^n = u(x_i,t^n)$. Consequently, we can write the finite Taylor series in the form
\begin{equation}
    \label{taylor}
u_i^n = U - \tau \partial_t U + \frac{\tau^2}{2} \partial_{tt} U + O(\tau^3) \,.
\end{equation}
Applying the standard Lax-Wendroff procedure
\begin{eqnarray}
\nonumber
\partial_t U = - \partial_x f(U) \,, \\[1ex]
\label{lw}
\partial_{tt} U = - \partial_{tx} f(U) = - \partial_{xt} f(U) = - \partial_{x} \left( f'(U) \partial_t U \right) =  \partial_{x} \left( f'(U) \partial_x f(U) \right) \,,
\end{eqnarray}
one obtains from (\ref{taylor}) that
\begin{equation}
    \label{lwf}
u_i^n = U + \tau \partial_{x} f(U) + \frac{\tau^2}{2}  \partial_x \left( f'(U) \partial_{x} f(U) \right)  + O(\tau^3) \,.
\end{equation}
One may obtain a numerical scheme from (\ref{lwf}) by approximating the spatial derivatives in (\ref{lwf}) using, e.g., finite differences. In such a way, a fully implicit scheme is obtained in the sense that it contains the value of $u$ at $t^n$ only on the left hand side of (\ref{lwf}). Such an approach can be found in the case of analogous explicit methods using more involved spatial discretizations in \cite{harten_uniformly_1997,shu_essentially_1998,qiu_finite_2003,leveque_finite_2004,toro_riemann_2009,seal_high-order_2014,tsai_two-derivative_2014}, and in the case of implicit method in, e.g., \cite{frolkovic2018semi,zeifang_two-derivative_2021}.

The first modification of the fully implicit approach (\ref{lwf}) to obtain a compact implicit scheme is to use only the first equality in (\ref{lw}) to obtain
\begin{equation}
\label{lwp}
u_i^n = U + \tau \partial_{x} f(U) - \frac{\tau^2}{2}  \partial_{tx} f(U)  + O(\tau^3) \,.
\end{equation}
Such an approach in even more general settings is used for analogous explicit schemes in \cite{tsai_two-derivative_2014, zorio_approximate_2017,li_two-stage_2019, carrillo2019compact, carrillo2021lax} and in the case of implicit schemes for ODEs in \cite{baeza_approximate_2020}.

To derive a second order finite difference scheme, one musts approximate the first occurrence of $\partial_x f$ in (\ref{lwp}) by a second order accurate scheme. We avoid a fixed stencil for such an approximation as we aim to derive a high resolution form to also approximate the discontinuous solutions of (\ref{cl}). Therefore, we introduce a convex combination of the second order central and upwind finite difference to approximate $\partial_x f$. We denote it by $\Delta^{\omega,2}_x$, where $\omega \in [0,1]$ is the free parameter, and 
\begin{equation}
    \label{delta2}
\partial_x f(u(x_i,t^*)) \approx \Delta^{\omega,2}_x f_i^* :=  (1-\omega) \frac{f_{i+1}^{*}-f_{i-1}^{*}}{2h} + \omega \frac{3 f_i^{*} - 4 f_{i-1}^{*} + f_{i-2}^{*}}{2h}
\end{equation}
with $*=n, n+1$, $f_i^* := f(u_i^*)$, and so on.

If the same spatial approximation is used to approximate $\partial_{tx} f$ together with the backward finite difference in time, one obtains a Crank-Nicolson type of the scheme 
$$
u_i^n = u_i^{n+1} + \tau \Delta^{\omega,2}_x f_i^{n+1} - \frac{\tau}{2}  \left( \Delta^{\omega,2}_x f_i^{n+1} - \Delta^{\omega,2}_x f_i^n \right)
$$
that can be written as 
$$
u_i^n = u_i^{n+1} +\frac{\tau}{2}  \left( \Delta^{\omega,2}_x f_i^{n+1} + \Delta^{\omega,2}_x f_i^n \right) \,.
$$
Clearly, such a scheme has the stencil of unknowns from $u_{i-2}^{n+1}$ up to $u_{i+1}^{n+1}$ if $\omega \in (0,1)$, see also \cite{coulette_implicit_2019,arbogast2021self,frolkovic2021report}. 

The final modification of the aforementioned approaches to obtain our compact implicit scheme is to use the following first order $\omega$-parametric approximation $\Delta_x^{\omega,1}$ of $\partial_x f$ instead of $\Delta_x^{\omega,2}$ to approximate $\partial_{tx} f$,
$$
\partial_x f(u(x_i,t^*)) \approx \Delta_x^{\omega,1} f_i^* := (1-\omega) \frac{f_{i+1}^* - f_i^*}{h} +  \omega \frac{f_i^* - f_{i-1}^*}{h} \,.
$$

Consequently, we obtain the following scheme
$$
u_i^n = u_i^{n+1} + \tau \Delta^{\omega,2}_x f_i^{n+1} - \frac{\tau}{2}  \left( \Delta^{\omega,1}_x f_i^{n+1} - \Delta^{\omega,1}_x f_i^n \right) 
$$
that can be simplified using
$$
\Delta_x^{\omega,2} f_i^{n+1} - \frac{1}{2} \Delta_x^{\omega,1} f_i^{n+1} = f_i^{n+1} - f_{i-1}^{n+1} - \frac{1}{2} \Delta_x^{\omega,1} f_{i-1}^{n+1}
$$
to the following final form,
\begin{eqnarray}
\label{final}
u_i^{n+1} + \frac{\tau}{h} \left(f_i^{n+1}-f_{i-1}^{n+1} -  \frac{1}{2} (\Delta_x^{\omega,1}f_{i-1}^{n+1}  - \Delta_x^{\omega,1} f_i^n) \right) = u_i^n \,.
\end{eqnarray}
%Note that () can be viewed as a linear combination of two 2nd order accurate approximations - firstly,  the central finite differences at $(x_{i-1/2},t^{n+1})$ and $(x_{i+1/2},t^n)$, and, secondly, the backward finite differences at $(x_{i+1/2},t^{n+1})$ and the central one at $(x_{i-1/2},t^n)$.
The advantage of this compact implicit second order accurate scheme is that it contains only the unknowns from $u_{i-2}^{n+1}$ up to $u_i^{n+1}$ for any fixed $\omega \in [0,1]$. In the case of linear advection equation, the scheme is unconditionally stable \cite{frolkovic_semi-implicit_2021}. 

The conservative form (\ref{step0}) can be obtained from (\ref{final}) by defining the numerical fluxes
\begin{equation}
    \label{flux0}
F_{i+1/2} = f_i^{n+1} - \frac{1}{2} \left( (1-\omega) (f_i^{n+1} - f_{i+1}^n) + \omega (f_{i-1}^{n+1} - f_{i}^n)\right) \,.
\end{equation}
It has clearly a form of the first order accurate flux (\ref{F}) corrected by the term
\begin{equation}
    \label{corr}
\frac{1}{2} \left( (1-\omega) (f_i^{n+1} - f_{i+1}^n) + \omega (f_{i-1}^{n+1} - f_{i}^n)\right) \,.
\end{equation}
Two special choices of $\omega$ in (\ref{corr}) should be mentioned. The case $\omega=1$ gives the fully upwinded form using only the available values of the numerical solution if the fast sweeping method is used. The second case with $\omega=0$ defines the correction with a ``central'' stencil that depends on the unknown value $u_i^{n+1}$.

Next, we propose a high resolution extension of the first order accurate numerical fluxes in (\ref{F}) based on (\ref{flux0}). We choose two tools in our high resolution compact implicit scheme to control the correction (\ref{corr}). The first choice is to define the solution-dependent values of $\omega$ for each numerical flux in the spirit of (Weighted) Essentially Non-Oscillatory schemes \cite{shu_essentially_1998}. The second choice is to decrease the update by (\ref{corr}) in (\ref{flux0}) by multiplying it with a factor having the value less than one.

To implement the two choices discussed before, the numerical flux functions in our high resolution compact implicit method will take the following parametric form,
\begin{eqnarray}
\label{F2a}
F_{i+1/2}^{+} = f_{i}^{+,n+1} - %\\[1ex]\nonumber
\frac{l_{i}}{2} \left( (1-\omega_{i}) (f_{i}^{+,n+1} - f_{i+1}^{+,n}) 
+ \omega_{i} (f_{i-1}^{+,n+1} - f_{i}^{+,n}) \right) \,, \\[1ex]
\label{F2b}
F_{i-1/2}^{-} = f_{i}^{-,n+1} - 
\frac{l_{i}}{2}\left( (1-\omega_{i}) (f_{i}^{-,n+1} - f_{i-1}^{-,n})
+ \omega_{i} ( f_{i+1}^{-,n+1} - f_{i}^{-,n}) \right)\,,
\end{eqnarray}
where the parameters $\omega_i \in [0,1]$ and $l_i \in [0,1]$ shall be chosen. These parameters are different in (\ref{F2a}) and (\ref{F2b}) \PF{for the same index $i$} (in fact, also in each time step), which we do not emphasize in the notation. 
%%

% The terms in (\ref{F2a}) and (\ref{F2b}) multiplied by $l_i$ can be seen as a second order correction of the first order accurate fluxes (\ref{F}) that are recovered by choosing $l_i=0$. We aim to use $l_i=1$ whenever possible, but in some cases we need to use $l_i \in [0,1)$, so the parameter $l_i$ plays a role of flux limiter \cite{}, see later.

% In the case of linear advection equation with constant parameters, the scheme is 2nd order accurate and unconditionally stable having no restriction on the choice of $\tau$ due to the stability, see the proofs in \cite{frolkovic_semi-implicit_2021}.
% We show next that for a fixed value of $\omega_i \equiv \bar \omega$ and $l_i \equiv 1$ the method (\ref{step0}) with either (\ref{F2a}) or (\ref{F2b}) is second order accurate for smooth solutions of (\ref{cl}). {\it stability for nonlinear case}

Note that the replacement of the first order accurate numerical fluxes (\ref{F}) in (\ref{step0}) by the high resolution forms (\ref{F2a}) and (\ref{F2b}) results again in a fully upwinded form in the implicit part of the scheme (\ref{step1}) - (\ref{step2}) for all particular choices of the parameters. Consequently, only the left hand sides of (\ref{step1}) and (\ref{step2}) contain the single unknown value $u_i^{n+1}$ if computed in the prescribed order of the fast sweeping method. The boundary conditions can be treated analogously to the first order scheme if we take into account that for $\omega_i=0$ one obtains an identical stencil in the implicit part. If the values $u_{-1}^n$ or $u_{I+1}^n$ are required by the high resolution numerical fluxes for the boundary nodes, they can be obtained by extrapolation.

\section{High resolution scheme}
\label{sec2}

In what follows, we propose a dependency of $\omega_i$ and $l_i$ in (\ref{F2a}) and (\ref{F2b}) on $u_i^{n+1}$ to obtain a Total Variation Diminishing (TVD) scheme, i.e., for $n\ge 0$ we keep the property (with some proper boundary conditions)
\begin{equation}
    \label{tvd}
\sum|u_i^n-u_{i-1}^n| \ge \sum  |u_i^{n+1}-u_{i-1}^{n+1}| \,.
\end{equation}
Note that the TVD schemes are treated extensively in literature, see, e.g.,  for explicit methods \cite{sweby1984high,harten1997high,gottlieb1998total,leveque_finite_2004,toro_riemann_2009,titarev2005weno} and for (explicit-)implicit methods \cite{harten_uniformly_1997,duraisamy_implicit_2007, dimarco2018second, michel2022tvd, puppo_quinpi_2022}.

For clarity of presentation, we derive the nonlinear TVD numerical scheme for the simplest case of the linear advection equation with a constant positive velocity $\bar v$, when $f(u) \equiv f^+(u)= \bar v u$. The Courant number is denoted by
\begin{equation}
    C = \frac{\bar v \tau}{h} \,. 
\nonumber
\end{equation}
In this case, the fluxes $F_{i\pm 1/2}$ in (\ref{step0}) take the form
\begin{eqnarray}
    \label{flux0adv}
F_{i+1/2} = \bar v \left(u_i^{n+1} - \frac{l_i}{2} \left( (1-\omega_i) (u_i^{n+1} - u_{i+1}^n) + \omega_i (u_{i-1}^{n+1} - u_{i}^n)\right) \right) \,, \\[1ex]
    \label{flux1adv}
F_{i-1/2} = \bar v \left(u_{i-1}^{n+1} - \frac{l_{i-1}}{2} \left( (1-\omega_{i-1}) (u_{i-1}^{n+1} - u_{i}^n) + \omega_{i-1} (u_{i-2}^{n+1} - u_{i-1}^n)\right) \right) \,.
\end{eqnarray}

\begin{rmk}
\label{rmktvd}
In what follows, we show that the scheme (\ref{step0}) with (\ref{flux0adv}) - (\ref{flux1adv}) can be formally expressed in the form
\begin{equation}
    \label{convex}
   u_i^{n+1} + c_{i-1} (u_i^{n+1} - u_{i-1}^{n+1})   = u_i^n \,,
\end{equation}
where the parameters $c_{i-1}$ depend on the numerical solution and on the parameters of the scheme. Afterwards, we define the values of $\omega_i$ and $l_i$ in such a way that $c_{i-1} \ge 0$. Consequently \cite{harten1997high,sweby1984high,puppo_quinpi_2022}, one can show that the scheme (\ref{convex}) is TVD. To do this, we sum the absolute values of the differences of (\ref{convex}) for $i$ and $i-1$ and use the reverse triangle inequality,
\begin{eqnarray}
    \nonumber
\sum_{i=2}^I |u_i^n-u_{i-1}^n| = \sum_{i=2}^I \left|(1+c_{i-1}) (u_i^{n+1}-u_{i-1}^{n+1}) - c_{i-2} ( u_{i-1}^{n+1} - u_{i-2}^{n+1})\right| \ge \\[1ex] \nonumber
 \sum_{i=2}^I \left( (1+c_{i-1}) |u_i^{n+1}-u_{i-1}^{n+1}| - c_{i-2} |u_{i-1}^{n+1} - u_{i-2}^{n+1}| \right) = \\[1ex] \nonumber -c_{0} |u_1^{n+1}-u_{0}^{n+1}| +  c_{I-1} |u_I^{n+1}-u_{I-1}^{n+1}| + 
\sum_{i=2}^{I} |u_i^{n+1}-u_{i-1}^{n+1}| .
\end{eqnarray}
Using, e.g., an assumption of the compact support of numerical solutions \cite{harten1997high,puppo_quinpi_2022}, i.e., $u_1^{n+1}=u_{0}^{n+1}$ and $u_I^{n+1}=u_{I-1}^{n+1}$, one obtains the TVD property (\ref{tvd}).
\end{rmk}

Let us transform (\ref{step0}) with (\ref{flux0adv}) - (\ref{flux1adv}) to the form (\ref{convex}). Suppose that $u_{i-1}^{n+1} \neq u_{i}^n$, then we can express the fluxes (\ref{flux0adv}) and (\ref{flux1adv}) in the form
\begin{eqnarray}
\label{upsip}
F_{i+1/2} = \bar v \left( u_{i}^{n+1} -\frac{l_i}{2}  \left( (1-\omega_{i}) \frac{u_{i}^{n+1} - u_{i+1}^{n}}{u_{i-1}^{n+1} - u_{i}^{n}}  + \omega_{i} \right)  (u_{i-1}^{n+1} - u_{i}^{n}) \right)\,,\\[1ex]
\label{upsim}
F_{i-1/2} = \bar v \left( u_{i-1}^{n+1} -\frac{l_{i-1}}{2} \left( 1-\omega_{i-1} + \omega_{i-1} \frac{u_{i-2}^{n+1} - u_{i-1}^{n}}{u_{i-1}^{n+1} - u_{i}^{n}} \right) (u_{i-1}^{n+1} - u_{i}^{n}) \right).
\end{eqnarray}

If we denote
\begin{equation}
    \label{r}
    r_{i-1} = \frac{u_{i-2}^{n+1} - u_{i-1}^n}{u_{i-1}^{n+1} - u_{i}^n}  \quad \Rightarrow \quad \frac{1}{r_i} = \frac{u_{i}^{n+1} - u_{i+1}^n}{u_{i-1}^{n+1} - u_{i}^n}
\end{equation}
and
\begin{equation}
    \label{psid}
    \Psi_{i-1} = 1-\omega_{i-1} + \omega_{i-1} r_{i-1} \quad\Rightarrow\quad
    \frac{\Psi_i}{r_i} = \frac{1-\omega_i}{r_i} + \omega_i \,,
\end{equation}
one can write the fluxes (\ref{upsip}) and (\ref{upsim}) in the concise form
\begin{eqnarray}
\label{upsiplim}
F_{i+1/2} = \bar v \left( u_{i}^{n+1} -\frac{l_i}{2} \frac{\Psi_i}{r_i}  (u_{i-1}^{n+1} - u_{i}^{n}) \right) \,, \\[1ex]
\label{upsimlim}
F_{i-1/2} = \bar v \left( u_{i-1}^{n+1} -\frac{l_{i-1}}{2} \Psi_{i-1} (u_{i-1}^{n+1} - u_{i}^{n}) \right)\,.
\end{eqnarray}

The values $\Psi_i$ can be viewed as the so-called flux limiters \cite{leveque_finite_2004,toro_riemann_2009,duraisamy_implicit_2007,puppo_quinpi_2022}. Although we aim to formulate our TVD scheme with a proper definition of the parameters $\omega_i \in [0,1]$ and $l_i \in [0,1]$, we have introduced the values $\Psi_i$  to make the following analysis and the comparison to previous publications easier.

The complete scheme (\ref{step0}) with (\ref{upsiplim}) and (\ref{upsimlim}) takes the form,
\begin{eqnarray}
    \label{steppsi}
     u_i^{n+1} - u_i^n + C \left( u_i^{n+1} - u_{i-1}^{n+1} \right)  -
     \frac{C}{2} \left(\frac{l_i \Psi_{i}}{r_{i}} - l_{i-1} \Psi_{i-1}\right) (u_{i-1}^{n+1} - u_{i}^n)  = 0 \,.
\end{eqnarray}
Now, using 
\begin{equation}
\nonumber
    \label{triv}
    u_{i-1}^{n+1} - u_{i}^n = (u_i^{n+1} - u_{i}^n) - (u_i^{n+1} - u_{i-1}^{n+1}) \,,
\end{equation}
the scheme (\ref{steppsi}) can be written in the form
% \begin{eqnarray}
% \label{tvdscheme}
% \left( 1 - \frac{C}{2} \left(\frac{l_i \Psi_{i}}{r_{i}} - l_{i-1} \Psi_{i-1}\right)\right)(u_i^{n+1} - u_i^n) + \\[1ex]\nonumber
% C \left( 1 + \frac{1}{2} \left(\frac{l_i \Psi_{i}}{r_{i}} - l_{i-1} \Psi_{i-1}\right) \right) (u_i^{n+1} - u_{i-1}^{n+1}) = 0 \,.
% \end{eqnarray}
\begin{eqnarray}
    \label{tvdscheme}
    u_i^{n+1} + \frac{C \left(
    1  + \frac{1}{2} \left(\frac{l_i \Psi_{i}}{r_{i}} - l_{i-1} \Psi_{i-1} \right)\right)}{1 - \frac{C}{2} \left(\frac{l_i \Psi_{i}}{r_{i}} - l_{i-1} \Psi_{i-1}\right)} \left(u_i^{n+1} - u_{i-1}^{n+1}\right) = u_i^n \,.
\end{eqnarray}
If the coefficients before $(u_i^{n+1} - u_{i-1}^{n+1})$ in (\ref{tvdscheme}) are well defined and nonnegative, then the scheme is TVD, see Remark \ref{rmktvd}. Note that if the denominator in (\ref{tvdscheme}) is zero for some $i$, the scheme (\ref{steppsi}) gives $u_i^{n+1}-u_{i-1}^{n+1}=0$, so the TVD property (\ref{tvd}) as shown in Remark \ref{rmktvd} is not destroyed.

In what follows, we formulate the inequalities for the values $\Psi_i$ and $l_i$ for which the assumptions of Remark \ref{rmktvd} are fulfilled for (\ref{tvdscheme}). We do it in two steps: first for $C\le 1$ and then for $C\ge 1$. 

In the first case where $C\le 1$, we can set $l_i=1$, so we do not use this kind of limiting. The values $\Psi_i$ can be defined using the TVD methods known for explicit methods \cite{leveque_finite_2004,kemm_comparative_2011,puppo_quinpi_2022}.
To have the positive coefficients in (\ref{tvdscheme}), it is enough to require 
\begin{eqnarray}
    \label{ineq1}
    -1 \le \Psi_{i-1} \le 2 \,, \\[1ex]
    \label{ineq2}
    \Psi_{i-1} - 2 \le \frac{\Psi_{i}}{r} \le \Psi_{i-1} + 2 \,,
\end{eqnarray}
where the inequalities in (\ref{ineq2}) must be satisfied for an arbitrary nonzero $r \in R$. Note that the inequalities for $\Psi_{i-1}$ in (\ref{ineq1}) are, in fact, required to fulfill (\ref{ineq2}) for two particular values of $\Psi_i$ that can occur: $\Psi_{i}=0$ and $\Psi_{i}=r_i$. 
%Note that for accuracy reasons, we require $\Psi(1)=1$ \cite{leveque_finite_2004,duraisamy_implicit_2007,puppo_quinpi_2022} {\it motivate or reformulate}.

There exist many variants of the construction of limiters that are typically given by defining a limiter function $\Psi=\Psi(r)$ such that $\Psi_i = \Psi(r_i)$. Analogously, using (\ref{psid}) one can define a function $\omega=\omega(r)$ such that $\omega_i=\omega(r_i)$.
One of the simplest choices is to define
\begin{equation}
\label{minmodomega}
\omega(r) = \left \{
\begin{array}{lr}
1 & |r| \le 1 \\[1ex]
0 & |r| > 1
\end{array}
\right . \,,
\end{equation}
that gives
\begin{equation}
\label{minmodpsi}
\Psi(r) = \left \{
\begin{array}{lr}
r & |r| \le 1 \\[1ex]
1 & |r| > 1
\end{array}
\right . \,.
\end{equation}
Clearly, for (\ref{minmodpsi}) the inequalities (\ref{ineq1}) and (\ref{ineq2}) are fulfilled.
In the case of fully explicit or fully implicit schemes, the choice (\ref{minmodomega}) can be viewed as the second order ENO reconstruction \cite{shu_essentially_1998,duraisamy_implicit_2007}. 

We propose the function $\omega=\omega(r)$ following a simple strategy of modified ENO schemes \cite{shu_numerical_1990}. That is, we choose a preferable constant value of $\omega$ to be used in (\ref{F2a}) that should be changed only if the TVD property is destroyed otherwise. In what follows, we choose the value $\omega=1$, i.e. $\Phi_i(r) = r$ due to  (\ref{psid}), which gives the most upwinded stencil of numerical fluxes in (\ref{F2a}) - (\ref{F2b}), see the discussion after (\ref{corr}). As shown for the linear advection equation in \cite{frolkovic2018semi,frolkovic_semi-implicit_2021}, this choice is preferable for large Courant numbers $C$ in (\ref{CN}) for accuracy reasons; see also the numerical experiments in Section \ref{ex01} for a comparison when solving the nonlinear Burgers' equation.

In particular, we define
\begin{equation}
    \label{omega1}
    \omega(r) = \left \{ \begin{array}{lr}
     \frac{1}{r-1}    &  2 \le r \\[1.5ex]
     \frac{2}{1-r}  & r \le -1 \\[1.5ex]
     1  & \hbox{otherwise}
    \end{array} \right.
\end{equation}
or, equivalently,
\begin{equation}
    \label{psi1}
    \Psi(r) = \left \{ \begin{array}{lr}
     2    &  2 \le r \\[1.5ex]
    -1  & r \le -1 \\[1.5ex]
     r  & \hbox{otherwise} \,.
    \end{array} \right.
\end{equation}
Clearly, if $\Psi_{i-1}$ and $\Psi_i$ are defined by (\ref{psi1}), then the inequalities (\ref{ineq1}) - (\ref{ineq2}) are fulfilled and the scheme is TVD. Note that if $\Psi(r) \neq r$ in (\ref{psi1}) then it takes the extremal values $-1$ and $2$ at the border of the TVD region.

% \begin{rmk}
% \label{rem-spec}
% To derive (\ref{tvdscheme}) we have supposed, among others, that $u_{i-1}^{n+1} \neq u_i^n$. As we show later, the case $u_{i-1}^{n+1} = u_i^n$ can happen only if $\omega_{i-1}=0$, when the scheme (\ref{step1}) takes the simpler form,
% \begin{eqnarray}
% \label{steppsispec1}
%  u_i^{n+1} - u_i^n + C \left( u_i^{n+1} - u_{i-1}^{n+1} - 
%  \frac{1}{2} (1-\omega_{i}) (u_i^{n+1}-u_{i+1}^n) \right) = 0 \,.
% \end{eqnarray}
% To fulfill the TVD property, we choose \PF{constant value $\omega_{i}=1$ in (\ref{steppsispec1}) which turns it into the first order scheme that is TVD}. If by chance the result of (\ref{steppsispec1}) is $u_i^{n+1}=u_{i+1}^n$, then $\omega_i$ in (\ref{steppsispec1}) can be arbitrarily chosen \PF{as it is multiplied by $0$, so we set $\omega_{i}=0$ to insure the property mentioned in the second sentence of this remark.}
% \end{rmk}

Next, we have to treat the case where $C\ge 1$. For that purpose, we use in (\ref{F2a}) and (\ref{F2b}) the factors $l_i \in [0,1]$ to limit the second order update (\ref{corr}) in numerical fluxes if necessary. To have positive coefficients in (\ref{tvdscheme}), we require more restrictive inequalities than (\ref{ineq1}) and (\ref{ineq2}), namely,
\begin{equation}
    \label{ineqe}
    -\frac{1}{C} \le l_{i-1} \Psi_{i-1} \le 2 \,, \quad -2 + l_{i-1}  \Psi_{i-1} \le \frac{l_{i} \Psi_{i}}{r} \le \frac{2}{C} + l_{i-1}  \Psi_{i-1} \,.
\end{equation}
Therefore, we define
\begin{equation}
    \label{omegae}
    \omega_{i} = \left \{ \begin{array}{lr}
     \frac{1}{r-1}    &  2 \le r \\[1.5ex]
     \frac{1+C}{C (1-r)}  & r \le -\frac{1}{C} \\[1.5ex]
     1  & \hbox{otherwise} 
    \end{array} \right.
\end{equation}
or, equivalently,
\begin{equation}
    \label{l}
     \Psi_{i}  = \left \{ \begin{array}{lr}
     2    &  2 \le r \\[1.5ex]
     -1/C  & r \le -\frac{1}{C} \\[1.5ex]
     r    &  \hbox{otherwise}.
    \end{array} \right. \,.
\end{equation}
Finally,
\begin{equation}
    \label{li}
l_i = \min \left \{ 1 , \max \left \{ 0, \frac{r_i}{\Psi_i}\left(\frac{2}{C}+l_{i-1}\Psi_{i-1}\right)\right \} \right \} \,.
\end{equation}
Using (\ref{omegae}) - (\ref{li}), one obtains the inequalities in (\ref{ineqe}) for arbitrary $C\ge 1$. Note that the definition of positive values for $\Psi_i$ in (\ref{l}) did not change with respect to (\ref{psi1}), so if $r_i \ge 0$ then the character of limiting using $\Psi_i$ is identical for any value of $C$.  

Finally, we comment on how to solve the algebraic equations (\ref{steppsi}) that are nonlinear due to the dependence of $r_i$ and $\Psi_i$ on the unknown value $u_i^{n+1}$. We propose it in the form of an iterative predictor-corrector procedure. We present all algorithmic details in Section \ref{sec3}, therefore, we mention here only some basic ideas. 

The scheme (\ref{step0}) is linear with the high resolution fluxes (\ref{F2a}) if some fixed values of parameters $\omega_i$ and $l_i$ are used. In all our numerical experiments, we compute a predicted value $u_i^{n+1,0} \approx u_i^{n+1}$ using (\ref{F2a}) or (\ref{F2b}) with $\omega_i=0$ and $l_i=1$. Once the predicted value $u_i^{n+1,0}$ is available, we compute the values of $r_i^0 \approx r_i$ using (\ref{r}) where we replace $u_i^{n+1}$ with $u_i^{n+1,0}$. Furthermore, we compute $\omega_i^0 \approx \omega_i$ using (\ref{omega1}) or (\ref{omegae}) by replacing $r_i$ with $r_i^0$. Similarly, the values $\Psi_i^0\approx \Psi_i$ from (\ref{psid}) and $l_i^0\approx l_i$ from (\ref{l}) are obtained; see Section \ref{sec3} for all details. Having the predicted values $\omega_i^0$ and $l_i^0$, we solve the linear algebraic equation (\ref{step0}) for $u_i^{n+1}$ where the numerical fluxes $F_{i+1/2}$ are defined with the fixed predicted values,
\begin{equation}
    \label{tvdflux}
    F_{i+1/2} = C \left(u_i^{n+1} - \frac{l_i^0}{2}\left((1-\omega_i^0) ( u_i^{n+1}-u_{i+1}^{n}) + \omega_i^0 (u_{i-1}^{n+1}-u_i^n)\right)\right) \,.
\end{equation}

In general, this procedure can be repeated using the computed corrected value as a new predictor. In our numerical experiments, we computed the corrector only once, and the numerical results for the chosen examples were satisfactory. In general, one cannot rely on such an approach, therefore, the algorithm in Section \ref{sec3} is formulated using, if necessary, several iterations of the predictor-corrector scheme using (\ref{tvdflux}). Note that the TVD property of our scheme is fulfilled only with the (nonlinear) high resolution fluxes in (\ref{F2a}) and (\ref{F2b}), therefore, if the corrected value is different from the predicted one, a small violation of the TVD property can occur in general.

Finally, we briefly comment on the treatment of the nonlinear flux function $f$. If, for example, $f'(u)\ge 0$ for $u \in R$, 
then the fluxes $F_{i\pm 1/2}$ in (\ref{step0}) take the form
\begin{eqnarray}
    \label{flux0n}
F_{i+1/2} = f_i^{n+1} - \frac{l_i}{2} \left( (1-\omega_i) (f_i^{n+1} - f_{i+1}^n) + \omega_i (f_{i-1}^{n+1} - f_{i}^n)\right) \,, \\[1ex]
    \label{flux1n}
F_{i-1/2} = f_{i-1}^{n+1} - \frac{l_{i-1}}{2} \left( (1-\omega_{i-1}) (f_{i-1}^{n+1} - f_{i}^n) + \omega_{i-1} (f_{i-2}^{n+1} - f_{i-1}^n)\right) \,,
\end{eqnarray}
where $f_i^n := f(u_i^n)$ and so on. The indicators $r_i$ are now determined by
\begin{equation}
    \label{rnl}
    r_{i} = \frac{f_{i-1}^{n+1} - f_{i}^n}{f_{i}^{n+1} - f_{i+1}^n} 
\end{equation}
giving for $f(u)=\bar v u$ the identical form as in (\ref{r}). Using (\ref{psid}), the numerical fluxes in (\ref{flux0n}) - (\ref{flux1n}) can be written again in the concise form
\begin{eqnarray}
\label{upsiplimnl}
F_{i+1/2} =  f_{i}^{n+1} -\frac{l_i}{2} \frac{\Psi_i}{r_i}  (f_{i-1}^{n+1} - f_{i}^{n}) \,, \\[1ex]
\label{upsimlimnl}
F_{i-1/2} =  f_{i-1}^{n+1} -\frac{l_{i-1}}{2} \Psi_{i-1} (f_{i-1}^{n+1} - f_{i}^{n}) \,.
\end{eqnarray}

Substituting (\ref{upsiplimnl}) - (\ref{upsimlimnl}) into (\ref{step0}) and using some straightforward algebraic manipulations, one obtains
\begin{eqnarray}
    \label{ltvdschemenl}
    u_i^{n+1} + 
    \frac{\frac{f_i^{n+1}-f_{i-1}^{n+1}}{u_i^{n+1}-u_{i-1}^{n+1}} \frac{\tau}{h} \left( 1 + \frac{1}{2} \left(\frac{l_{i}\Psi_{i}}{r_{i}} - l_{i-1}\Psi_{i-1}\right) \right)}{1 - \frac{1}{2} \frac{f_i^{n+1}-f_i^n}{u_i^{n+1}-u_i^{n}} \frac{\tau}{h} \left(\frac{l_{i}\Psi_{i}}{r_{i}} - l_{i-1}\Psi_{i-1}\right)} (u_i^{n+1} - u_{i-1}^{n+1}) = u_i^n \,.
\end{eqnarray}
The equation (\ref{ltvdschemenl}) is similar to (\ref{tvdscheme}) with one important difference that the constant Courant number $C$ in (\ref{tvdscheme}) is replaced in (\ref{ltvdschemenl}) by nonlinear terms. In theory, if some estimate of the maximal value of
$$
\frac{f_i^{n+1}-f_i^n}{u_i^{n+1}-u_i^{n}}
$$
is available, the definitions of parameters $\omega_i$ and $l_i$ can be motivated analogously to (\ref{tvdscheme}).
In the next section, we describe all the details of the method using the high resolution compact implicit scheme for the scalar (nonlinear) hyperbolic equation.

\section{The high resolution compact implicit scheme}
\label{sec3}

 For simplicity, we suppose that the solution $u$ of (\ref{cl}) has a compact support, so we can set
\begin{equation}
\nonumber
    \label{compactl}
    u_0^{n+1} = u_0^n \,, \,\, u_1^{n+1} = u_1^{n} \,,
\end{equation}
and 
\begin{equation}
    \label{compactr}
    \nonumber
    u_I^{n+1} = u_I^n \,, \,\, u_{I-1}^{n+1} = u_{I-1}^{n} \,.
\end{equation}
Furthermore, we assume that Courant numbers $C^+ \ge 0$ and $C^- \ge 0$ are available such that
\begin{equation}
    \label{CN}
    \nonumber
    \frac{\tau}{h} \max_u \frac{d}{du} f^+(u) \le C^+ \quad \hbox{and} \quad \frac{\tau}{h} \max_u \frac{d}{du} f^-(u) \ge -C^- \,.
\end{equation}
\PF{Let $\epsilon>0$ denote a chosen small enough constant for which some terms in the following algorithm, if smaller than $\epsilon$, are considered to be zero effectively}. The following default values are used in each time step, if not defined otherwise: $\omega_i=0$, $l_i=1$, and $\Psi_i=1$.
The scheme (\ref{step1}) using (\ref{F2a}) for $i=2,3,\ldots,I-2$ or (\ref{step2}) using (\ref{F2b}) for $i=I-2,I-3,\ldots,2$ is then iteratively solved at the $n$-th time step as follows.
\vspace{2ex}
\begin{enumerate}

\item Compute 
\begin{equation}
    \label{fupw}    \nonumber
    \Delta^{up} = f_{i-1}^{+,n+1} - f^{+,n}_i \quad \hbox{or} \quad
    \Delta^{up} = f_{i+1}^{-,n+1} - f_i^{-,n} \,.
\end{equation}

If $|\Delta^{up}|\le \epsilon$ then set $\omega_i=1$ and solve the algebraic equation (\ref{step1}) or (\ref{step2}) for the unknown $u_i^{n+1}$. Continue with the step 1 for $i+1$ or $i-1$.\\

\item If $|\Delta^{up}| > \epsilon$ then set an initial guess $\mathrm{u}^{0} \approx u^{n+1}_i$  by using, e.g., the first order accurate numerical fluxes (\ref{F}) in (\ref{step1}) or (\ref{step2}) and solving for $\mathrm{u}$ the equation
\begin{equation}
     \nonumber
    \mathrm{u} + \frac{\tau}{h} f^+(\mathrm{u}) = u_i^n + \frac{\tau}{h} F^{+}_{i-1/2}
\end{equation}
or
\begin{equation}
    \nonumber
    \mathrm{u} - \frac{\tau}{h} f^-(\mathrm{u}) = u_i^n - \frac{\tau}{h} F^{-}_{i+1/2} \,,
\end{equation}
or by using the second order accurate numerical fluxes with $\omega_i=0$ and $l_i=1$,
\begin{equation}
    \label{fo1} 
    \mathrm{u} + \frac{\tau}{2 h} f^+(\mathrm{u}) = u_i^n + \frac{\tau}{h} F^{+}_{i-1/2} - \frac{\tau}{2 h}f^+(u_{i+1}^n)
\end{equation}
and
\begin{equation}
    \label{fo2}
    \mathrm{u} - \frac{\tau}{2 h} f^-(\mathrm{u}) = u_i^n - \frac{\tau}{h} F^{-}_{i+1/2} + \frac{\tau}{2 h} f^-(u_{i-1}^n)  \,.
\end{equation}
\item For the value $\mathrm{u}^k \approx u_i^{n+1}$ for some $k\ge 0$ compute
\begin{equation}
    \label{fdw}\nonumber
    \Delta^{dw,k} =  f^+(\mathrm{u}^k) - f_{i+1}^{+,n} \quad \hbox{or} \quad
    \Delta^{dw,k} =  f^-(\mathrm{u}^k) - f_{i-1}^{-,n} \,.
\end{equation}
If $|\Delta^{dw,k}|\le \epsilon$ then proceed with the step 5.

\item If $|\Delta^{dw,k}| > \epsilon$ then compute
\begin{equation}
    \label{main}\nonumber
    r^k = \frac{\Delta^{up}}{\Delta^{dw,k}} \,, 
\end{equation}
and
\begin{eqnarray}
    \label{omeganew}
    \omega_i^k = \left \{ \begin{array}{lr}
     \frac{1}{r^k-1}    &  2 \le r^k \\[1.5ex]
     \frac{1+C}{C (1-r^k)}  & r^k \le -\frac{1}{C} \\[1.5ex]
     1  & \hbox{otherwise} ,
    \end{array} \right. 
\end{eqnarray}
with $C=\max\{1,C^+\}$ or $C=\max\{1,C^-\}$, respectively. Furthermore,
\begin{equation}
    \label{psiM}\nonumber
    \psi_i^k = 1 - \omega_i^k + \omega_i^k r^k \,
\end{equation}
and if $\psi_i^k \neq 0$ then
\begin{equation}
    \label{lnewl}
    l_i^k = \min \{ 1 , \max \{ 0, \frac{r^k}{\psi_i^k}\left(\frac{2}{C}+l_{i-1}\psi_{i-1}\right)\} \} \,,
\end{equation}
or
\begin{equation}
    \label{lnewr}\nonumber
    l_i^k = \min \{ 1 , \max \{ 0, \frac{r^k}{\psi_i^k}\left(\frac{2}{C}+l_{i+1}\psi_{i+1}\right)\} \}\,.
\end{equation}

\item Having established the $k$-th estimates of the parameters $\omega_i$ and $l_i$, we solve the algebraic equation (\ref{step1}) or (\ref{step2}) with (\ref{F2a}) or (\ref{F2b}) and denote its solution by $u_i^{n+1,k+1}$. If a chosen stopping criterion is fulfilled, \PF{e.g., $|u_i^{n+1,k+1} - u_i^{n+1,k}|<\epsilon$ or a prescribed number of corrector steps is reached}, we set $u_i^{n+1}=u_i^{n+1,k+1}$ and proceed with the step 1 for $i+1$ or $i-1$. If not, we proceed with step 3.\\[0ex]
\end{enumerate}

We note that to improve accuracy, one can replace $C^+$ in (\ref{omeganew}) by $C_i^k$ defined by
\begin{equation}
    \label{Cnew}\nonumber
    C_{i}^k = \left \{ \begin{array}{lr}
     \frac{f^+(u_i^{n+1,k})-f^+(u_i^{n})}{u_i^{n+1,k}-u_i^{n}}    &  u_i^{n+1,k}\neq u_i^{n} \\[2ex]
     \frac{d}{du}f^+(u_i^{n+1,k})  & u_i^{n+1,k} = u_i^{n}
    \end{array} \right.
\end{equation}
and analogously for $C^-$.

\section{Hyperbolic systems}
\label{sec4}

Concerning the systems of hyperbolic equations, one has to take the steps defined for the scalar case in the previous section for each component of the system. Similarly to experiences in several papers \cite{leveque_finite_2004,duraisamy_implicit_2007,carrillo2021lax, boscarino_high_2022}, we prefer to express the higher order update (\ref{corr}) of the first order numerical flux in (\ref{F}) with the help of characteristic variables and speeds (the eigenvalues).

We consider the system (\ref{cl}) for ${\bf f} : R^m \rightarrow R^m$,
\begin{equation}
    \label{sys}
    \partial_t {\bf u} +\partial_x {\bf f}({\bf u}) = 0 \,,
\end{equation}
where we suppose that the Jacobian matrix ${\bf f}'({\bf u})$ has only nonnegative real eigenvalues $\lambda^p$, $p=1,2,\ldots,m$. The systems with nonpositive eigenvalues are treated analogously, the general case is solved using the fractional step method as explained before in (\ref{step1}) - (\ref{step2}) \cite{lozano_implicit_2021}.

Let the columns of the matrix $R=R({\bf u})$ be given by the eigenvectors ${\bf r}^p$, $p=1,2,\ldots,m$. Due to the hyperbolicity of (\ref{sys}), the matrix $R$ is regular for each considered value of ${\bf u}$.
Let ${\bf u}$ be the last estimate (or predictor) of ${\bf u}^{n+1}_i$ and let $R^{-1}$ be the inverse matrix of $R({\bf u})$. We express the term ${\bf f}_i^{n+1,0} - {\bf f}_{i+1}^n$ and the term ${\bf f}_{i-1}^{n+1} - {\bf f}_{i}^n$ in the second order update (\ref{corr}) of (\ref{flux0}) as a linear combination of eigenvectors using
\begin{equation}
    \label{ev}\nonumber
   \boldsymbol{\alpha}_i = R^{-1} \cdot \left(  {\bf f}_i^{k,n+1} - {\bf f}_{i+1}^n  \right) \,, \quad
   \boldsymbol{\beta}_i = R^{-1} \cdot \left(  {\bf f}_{i-1}^{n+1} - {\bf f}_{i}^n  \right) \,.
\end{equation}
Subsequently, the high resolution fluxes in (\ref{F2a}) are defined by
\begin{eqnarray}
    \label{sisys2}
     {\bf F}_{i+1/2} = {\bf f}_i^{n+1} - \frac{1}{2}  
     \sum_p l_i^p \left( (1-w^p_i) \alpha^p_i  + w_i^p \beta^p_i \right) {\bf r}^p \,,
\end{eqnarray}
where the weights in ${\bf w}_i =(w^1_i,w^2_i,\ldots,w^m_i)$ and ${\bf l}_i=(l^1_i,l^2_i,\ldots,l^p_i)$ are now associated with the components in ${\boldsymbol \alpha}_i$ and ${\boldsymbol \beta}_i$. 
%\begin{eqnarray}
%    \label{sisys}
%     {\bf F}_{i+1/2}^{k,n+1} = {\bf f}_i^{n+1} - \frac{1}{2} \left( (1-w_i^1) \alpha_i^1 \lambda^1 r^1 + (1-w_i^2) \alpha_i^2 \lambda^2 r^2
%     + w_i^1 \alpha_{i-1}^1 \lambda^1 r^1 + w_i^2\alpha_{i-1}^2 \lambda^2 r^2 \right) \right. - \ldots
%\end{eqnarray}

Having the form (\ref{sisys2}), the high resolution approach of the scalar case is used for each component of the system with the indicators ${\bf r}_i$ defined by
\begin{equation}
    \label{rsys}
    r^p_i = \frac{\beta^p_i}{\alpha^p_i} \,, \,\, p=1,2,\ldots,m \,.
\end{equation}
Furthermore, the Courant numbers $C^+$ in (\ref{omeganew}) are replaced by the corresponding values of the eigenvalues $\lambda^p$ for each component.

% In particular, let ${\bf u}={\bf u}^{k,n+1}_i$ denote the last estimate of ${\bf u}^{n+1}_i$. Our aim is to suggest the values of ${\bf \omega}_i$ such that
% \begin{equation}
%     \label{ineqsys}
%     -\frac{1}{\lambda^p} < \Psi_{i}^p \le 2 \,, \quad -2 + l_{i-1}  \Psi_{i-1}^p \le \frac{l_{i} \Psi_{i}^p}{r^p} \le \frac{2}{\lambda^p} + l_{i-1}^p  \Psi_{i-1}^p \,.
% \end{equation}

% To do so we define
% \begin{equation}
%     \label{rsys}
%     r^p = \frac{\beta^p}{\alpha^p} \,, \,\, p=1,2,\ldots,m
% \end{equation}
% and
% \begin{equation}
%     \label{omegasys}
%     \omega_{i}^p = \left \{ \begin{array}{lr}
%      \frac{1}{r^p-1}    &  2 \le r^p \\[1.5ex]
%      \frac{2+\lambda^p}{\lambda^p (1-r^p)}  & r^p \le -\frac{2}{\lambda^p} \\[1.5ex]
%      1  & \hbox{otherwise} .
%     \end{array} \right.
% \end{equation}
% and
% \begin{equation}
%     \label{psisys}
%      \Psi_{i}^p  = 1-\omega^p_i + \omega^p_i r^p
% \end{equation}
% and
% \begin{equation}
%     \label{lsys}
%     l_{i}^p = \min \left \{ 1 , \max\left \{0,\frac{r^p}{\Phi_i^p} \left(\frac{2}{\lambda^p}+l_{i-1}^p\Psi_{i-1}^p\right)\right \} \right \} \,.
% \end{equation}

% As noted before, the coefficients $\alpha$ and $\beta$ are depending on $u_i^{n+1}$. An initial guess for $u_i^{n+1}$ can be obtained by choosing $w_i^p=1$ in (\ref{sisys}).

\section{Numerical experiments}
\label{sec5}

In what follows, we illustrate numerical resolutions of the proposed implicit scheme for several standard test problems taken from the literature. The implementation is carried out using Mathematica software \cite{Mathematica}. If available, we present exact solutions to the examples that are then used to set the boundary conditions.

All examples are computed with the numerical scheme having the form of fractional step method (\ref{step1}) - (\ref{step2}). Note that if all characteristic speeds are nonnegative, it is enough to compute only the first step (\ref{step1}). The numerical fluxes of the first order accurate scheme are given by (\ref{F}), and the numerical fluxes of the high resolution scheme are defined in (\ref{F2a}) - (\ref{F2b}). The predictor is always computed using (\ref{fo1}) - (\ref{fo2}), the parameters $\omega_i$ and $l_i$ are then computed by (\ref{omeganew}) and (\ref{lnewl}) - (\ref{lnewr}), respectively. Only one corrector step is used by solving (\ref{step1}) - (\ref{step2}) with the predicted values $\omega_i$ and $l_i$ in (\ref{F2a}) - (\ref{F2b}).

When computing the examples for Burgers' equation with $f(u)=u^2/2$, we use the approach of \cite{lozano_implicit_2021} when the splitting (\ref{split}) is obtained by  
\begin{equation}
    \label{bsplit}\nonumber
    f^+(u):=\frac{1}{2}\left( f(u)+\left|u\right|\frac{u}{2}\right) \,, \quad 
    f^-(u):=\frac{1}{2}\left( f(u)-\left|u\right|\frac{u}{2}\right) \,.
\end{equation}

\subsection{Linear advection}
\label{ex03}

To illustrate the TVD property of our scheme, we solve the test example \cite{jiang_efficient_1996,balsara2000monotonicity,borges_improved_2008,kemm_comparative_2011} with non-smooth solutions for the advection with constant unity speed. This example is often used to judge the behavior of non-oscillatory numerical schemes in the literature. The computational domain is the interval $[-1,1]$. The initial condition consists of four different segments - a Gaussian, a triangle, a square-wave and a semi-ellipse, and it takes the following form \cite{balsara2000monotonicity}
\begin{equation}
    \label{initadv}
    u(x,0) = \left \{
\begin{array}{lr}
 \frac{1}{6}\left(G(x,\beta,z-\delta) + G(x,\beta,z+\delta) + 4 G(x,\beta,z) \right)      &-0.8 \leq x \leq -0.6\\[1ex]
 1                          &-0.4 \leq x \leq -0.2\\[1ex]
 1-|10(x-0.1)|              & 0.0 \leq x \leq 0.2\\[1ex]
 \frac{1}{6}\left(F(x,\alpha,a-\delta) + F(x,\alpha,a+\delta) + 4 F(x,\alpha,a) \right)     & 0.4 \leq x \leq 0.6\\[1ex]
0   & \hbox{otherwise} \,,
\end{array}
    \right .
\end{equation}
where
\begin{equation}
\label{adv_g}
    G(x,\beta,z)=e^{-\beta(x-z)^2} \, , \,\,\, 
    F(x,\alpha,a)=\sqrt{\max(1-\alpha^2(x-a)^2,0)} \, .
\end{equation}
The constants are given by
\begin{equation}\nonumber
    a=0.5, \quad z=-0.7, \quad \delta=0.005, \quad \alpha=10, \quad \beta=\frac{\log 2}{36\delta^2} \, .
\end{equation}

\begin{figure}
    \centering
    \includegraphics[width=6.0cm]{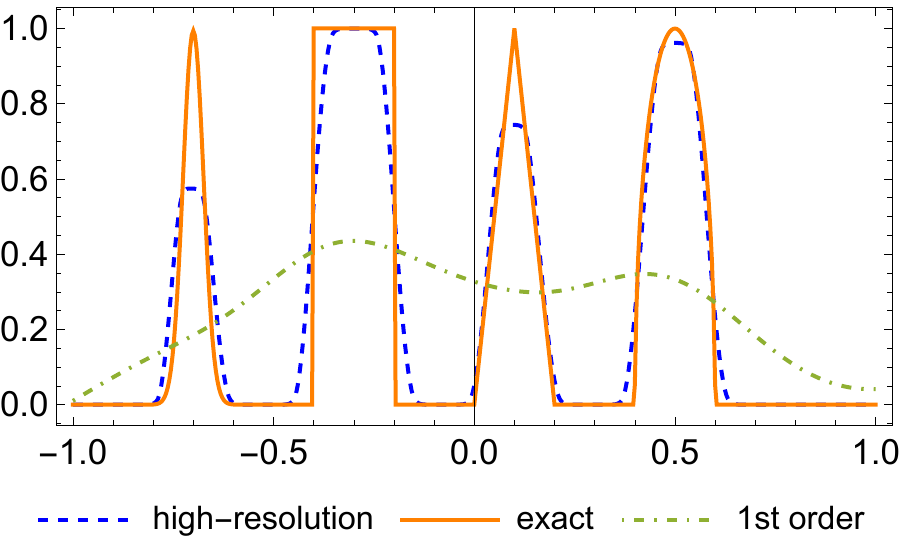}
    \includegraphics[width=6.0cm]{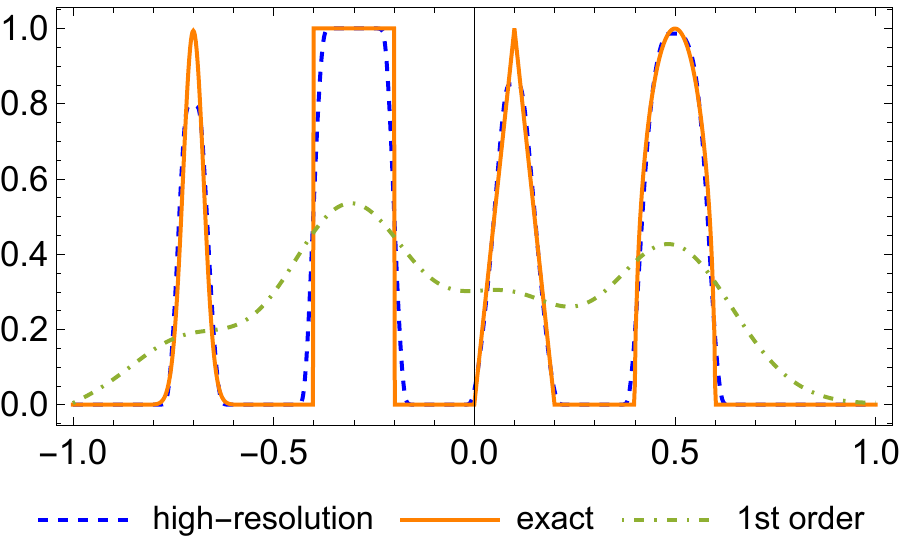}
    \caption{The comparison of the exact and the numerical solutions obtained with the first order scheme and the high resolution scheme for the example in Section \ref{ex03}. The left picture is obtained for $I=500$ and the right one for $I=1000$ after $125$ and $250$ time steps, respectively. The Courant number is always $4$, so $\tau = 4 h$.}
    \label{fig:lin}
\end{figure}

The problem is solved with $C=4$, so $\tau=4h$, and the numerical solutions are shifted backward after each time step to return to the initial position.
For a visual evaluation of the results obtained with the high resolution method see Figure \ref{fig:lin}, where a clear improvement with respect to the first order scheme can be seen. Moreover, no over- or undershootings with magnitudes larger than rounding errors are observed. Such results are obtained with time steps larger than typically allowed for explicit schemes.

\subsection{The smooth solution of Burgers' equation}
\label{ex01}

\begin{figure}
    \centering
    \includegraphics[width=6.0cm]{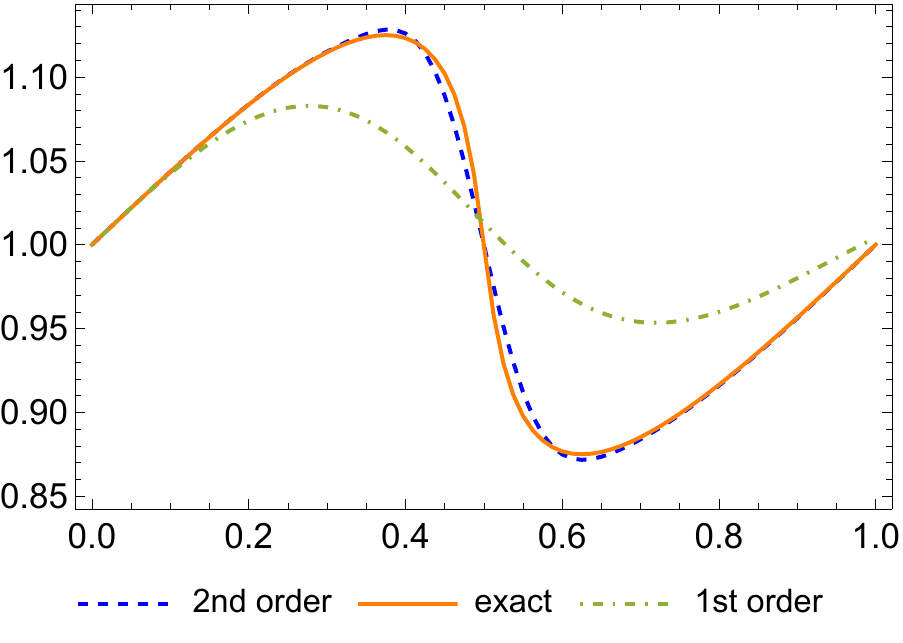}
    \includegraphics[width=6.0cm]{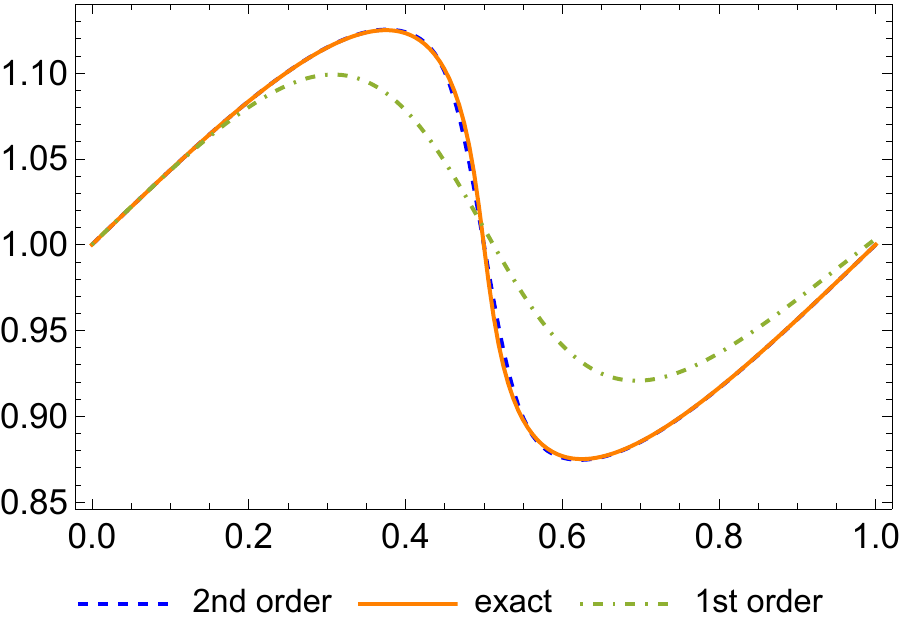}
    \caption{The comparison of the exact and the numerical solutions obtained with the first order scheme and the second order scheme with $\omega = 1$ for the example in Section \ref{ex01}. The left picture is obtained for $I=80$ and the right one for $I=160$ after $20$ and $40$ time steps, respectively. The maximal Courant number is always $4.5$ with $\tau=4 h$.}
    \label{fig:nsin}
\end{figure}

\begin{table}[ht]
    \small
	\begin{center}
	\begin{tabular}{c c c c c c c c c c}
	\hline
 $I$ & $1st$ order & EOC & $\omega=0$ & EOC & $\omega=1/2$ & EOC & $\omega=1$ & EOC\\ 
	\hline
 40 & 0.04214 & - & 0.01357 & - & 0.00761 & - & 0.00342 & - \\
 80 & 0.02525 & 0.74 & 0.00428 & 1.66 & 0.00230 & 1.73 & 0.00091 & 1.91 \\ 
 160 & 0.01419 & 0.83 & 0.00121 & 1.81 & 0.00064 & 1.84 & 0.00021 & 2.08 \\
 320 & 0.00768 & 0.89 & 0.00033 & 1.89 & 0.00017 & 1.92 & 0.00005 & 2.17 \\
	\hline 
	\end{tabular}
	\end{center}
	\caption{The errors and the experimental order of convergence (EOC) for the first order scheme and the second order scheme with three different values of parameter $\omega=0, 0.5, 1$ for the example in Section \ref{ex01}. The maximal Courant number is $4.5$ .}
	\label{tabsin}
\end{table}

In this example, we test the scheme (\ref{step0}) with the second order accurate numerical fluxes (\ref{flux0}) for the fixed values of $\omega$ in the case of a smooth solution of Burgers' equation. That is, we set
\begin{equation}
    \label{sininit}\nonumber
    f(u) = \frac{u^2}{2} \,, \quad u(x,0) = 1+\frac{1}{8} \sin(2\pi x) \,, \,\, x \in [0,1] \,,
\end{equation}
and we solve the equation for $t \in [0,1]$. The exact solution is computed numerically using the method of characteristics by solving the algebraic equations for the unknowns $u=u(x_i,t^n)$
$$
 u = 1 + \frac{1}{8} \sin (2\pi (x_i-u t^n)) \,.
$$

In Figure \ref{fig:nsin}, a comparison of the exact and numerical solutions obtained with the first and the second order method is given at the final time for two grids with $I=80$ and $I=160$ and $\tau = 4 h$. 
The global $l_1$ discrete error in time and space
\begin{equation}
    \label{l1}
    E_I^N := h \tau \sum \limits_{i=0}^I \sum \limits_{n=1}^N |u_i^n - u(x_i,t^n)| 
\end{equation}
is presented for all cases in Table \ref{tabsin}. The expected EOC for the first and the second order schemes is confirmed with the most accurate results delivered by the second order scheme using $\omega=1$ for the chosen maximal Courant number $4.5$.

\subsection{Slowly moving shock of Burgers' equation}
\label{ex02}

\begin{figure}
    \centering
    \includegraphics[width=6.0cm]{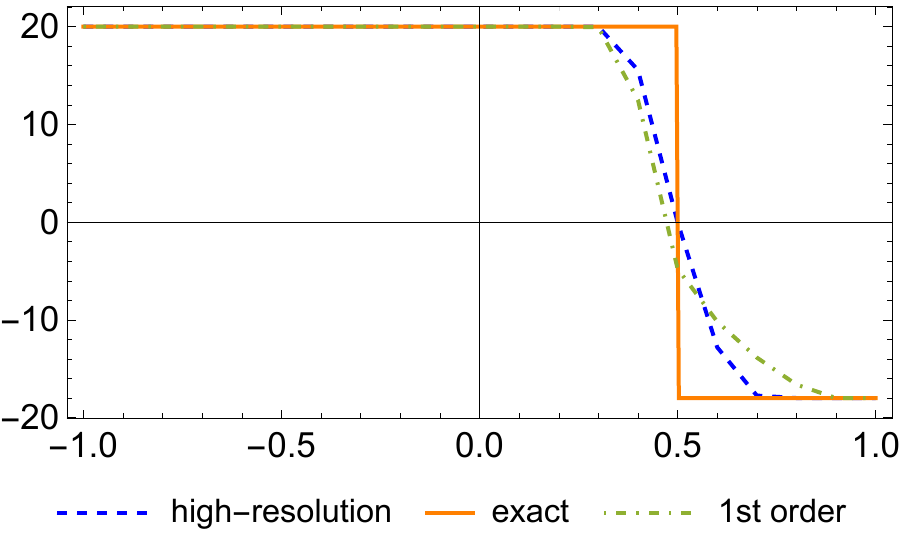}
    \includegraphics[width=6.0cm]{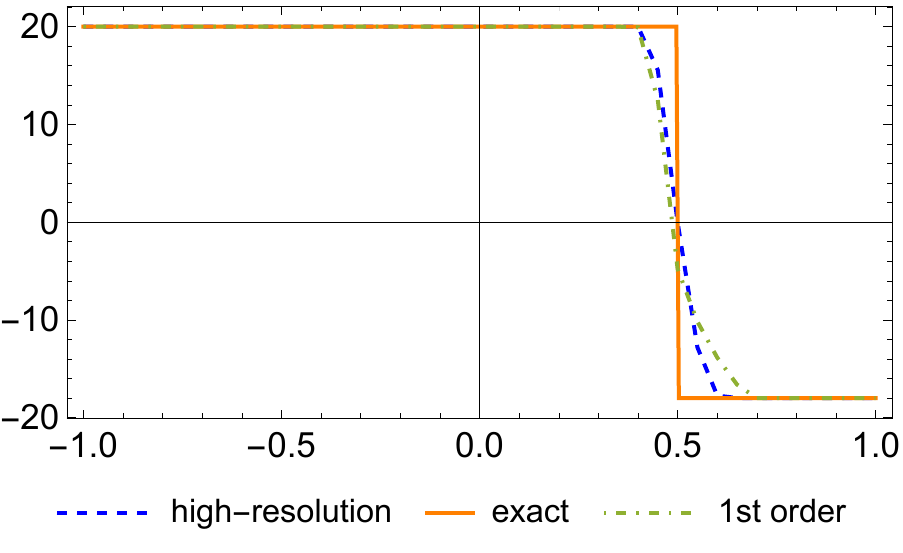}
     \caption{The comparison of the exact and the numerical solutions obtained with the first order method and the high resolution method for the example in Section \ref{ex02}. The left picture is obtained for $I=20$ and the right one for $I=40$ after $40$ and $80$ time steps, respectively. The maximal Courant number is always $10$.}
    \label{fig:slow}
\end{figure}

Inspired by \cite{lozano_implicit_2021}, we present numerical solutions obtained with the high resolution scheme for the Riemann problem with a slowly moving shock. This example can be seen as the simplest illustration when the implicit scheme can be computationally more efficient than some analogous explicit schemes.

The initial discontinuity of the piecewise constant function is placed at $x=-0.5$ with the left value $u_L=20$ and the right value $u_R=-18$. Consequently, the shock speed is equal to $1$ \cite{leveque_finite_2004}, so the piecewise constant profile is preserved with the moving position $x=-0.5+t$ of the discontinuity. We present the comparison of the first order and the high resolution scheme in $t=1$ in Figure \ref{fig:slow} for two rather coarse meshes with $I=20$ and $I=40$. We use the time step $\tau=h/2$ that is much larger than the explicit schemes would typically allow, because it corresponds to the maximal Courant number equal to $10$. One can observe on the coarse grid a significantly improved approximation of the shock speed for the numerical solution obtained with the high resolution scheme when compared with the first order scheme.

\subsection{Burgers' equation with interacting shock and rarefaction}
\label{excomplex}

\begin{figure}
    \centering
    \includegraphics[width=6.0cm]{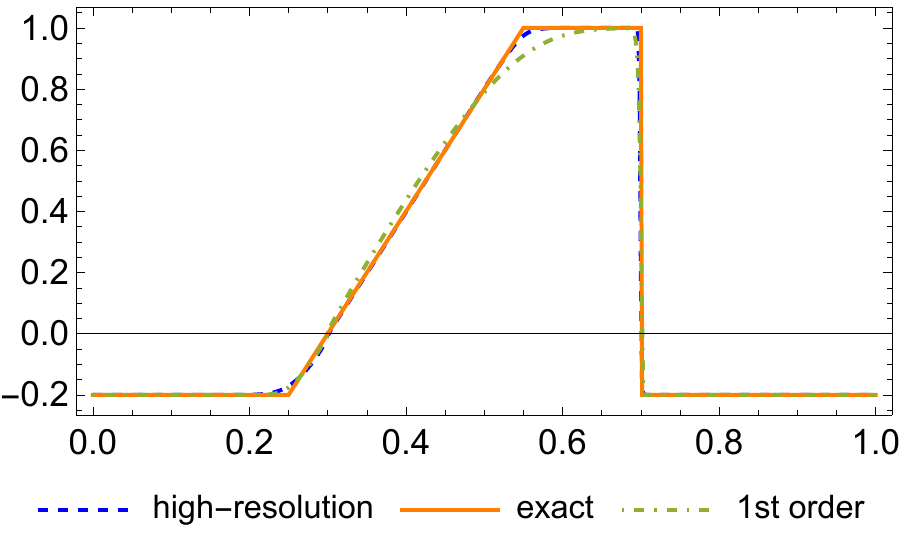}
    \includegraphics[width=6.0cm]{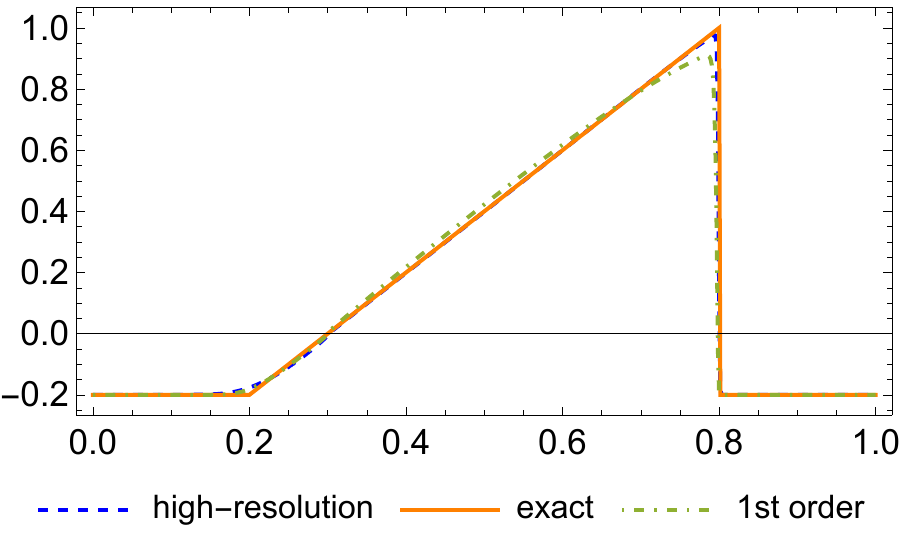}\vspace{.3cm}

    \includegraphics[width=6.0cm]{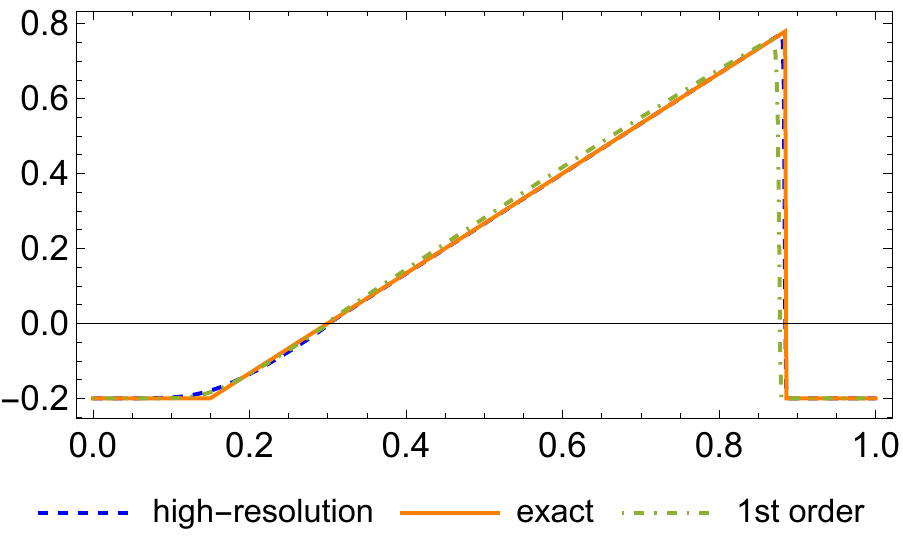}
    \includegraphics[width=6.0cm]{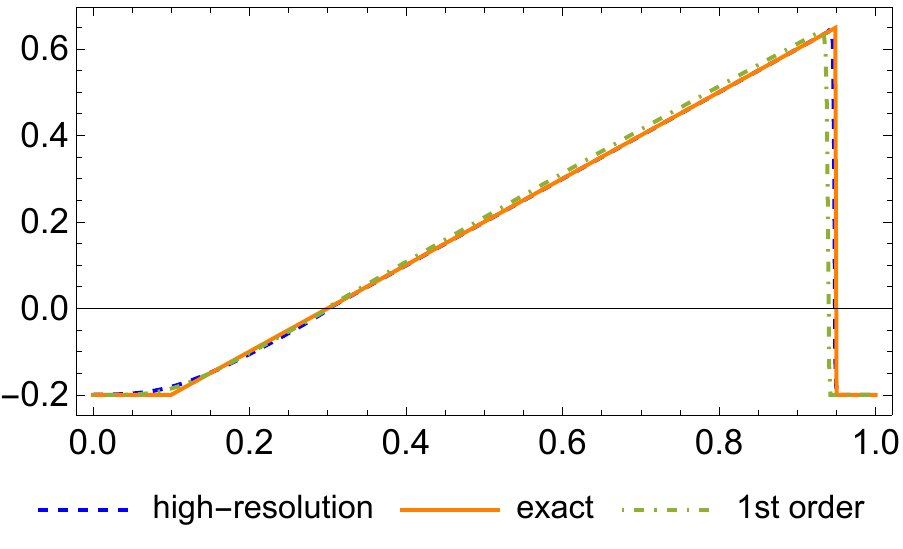}
    
    \caption{The comparison of the exact solution with the numerical solutions obtained with the first order scheme and the high resolution scheme for the example in Section \ref{excomplex}. The first row presents the results at $t=0.25$ (left) and $t=0.5$ (right) and the second row at $t=0.75$ and $t=1$. The number of mesh points is $I=640$, the maximal Courant number is $4$ with $\tau =4 h$.}
    \label{fig:complex}
\end{figure}

The last example of Burgers' equation is taken from \cite{lozano_implicit_2021} and contains all typical features of the solutions of Riemann problems for this equation. The purpose is to test the behavior of the high resolution scheme for the approximations of a nontrivial interaction between a shock and a rarefaction wave when choosing time steps larger than typically allowed for explicit schemes.

The initial condition is given by
\begin{equation}
    \label{initcomplex}
    \nonumber
    u(x,0) = \left \{
\begin{array}{lr}
 1    &  0.3 < x < 0.6\\[1ex]
-0.2   & \hbox{otherwise}
\end{array}
    \right . 
\end{equation}
and the exact solution is taken from \cite{lozano_implicit_2021} which is defined for $t \in (0,0.5)$ by:
\begin{equation}
    \label{complexex}
    u(x,t) = \left \{\begin{array}{lr}
       \frac{x-0.3}{t}  &  0.3 - 0.2 t \le x \le 0.3 + t \\[1ex]
       1 & 0.3 + t \le x < 0.6 + 0.4 t \\[1ex]
       -0.2  & \hbox{otherwise}.
    \end{array}
    \right.
\end{equation}

At time $t=0.5$ the end points of the rarefaction wave and the shock wave merge, and the solution evolves further with a triangular profile for $t\ge 0.5$
\begin{equation}
    \label{complexex}
    u(x,t) = \left \{\begin{array}{lr}
       \frac{x-0.3}{t}  &  0.3 - 0.2 t \le x < 0.3 - 0.2 t + 0.6 \sqrt{2 t} \\[1ex]
       -0.2  & \hbox{otherwise}.
    \end{array}
    \right.
\end{equation}

One can see that the high resolution method gives satisfactory results for this complex example even if the maximal Courant number is $4$. In Figure \ref{fig:complex} one can see that the first order scheme approximates the exact solution at $t=0.5$ with a visibly larger error than the high resolution scheme. Probably, this imprecision is a reason why the position of the shock moving with variable speed for $t>0.5$ is significantly better approximated with the high resolution method. Both methods converge to the exact solution with respect to the error defined in (\ref{l1}), see Table \ref{tab}.

\begin{table}[ht]
	\begin{center}
	\begin{tabular}{ c c l l l l }
	\hline
 $I$ & $N$ & $E_I^N$ & EOC & $E_I^N$ & EOC \\ 
	\hline
 160 & 40 & 0.01042 & - & 0.0374 & - \\
 320 & 80 & 0.00564 & 0.85 & 0.0235 & 0.67 \\
 640 & 160 & 0.00314 & 0.84 & 0.0144 & 0.71 \\
 1280 & 320 & 0.00175 & 0.84 & 0.0087 & 0.73 \\
	\hline 
	\end{tabular}
	\end{center}
	\caption{The numerical errors (\ref{l1}) with EOC for the example \ref{excomplex}. The third and the fourth columns are for the high resolution method, the fifth and the sixth ones for the first order method. The maximal Courant number is $4$ with $\tau =4 h$.}
	\label{tab}
\end{table}

\subsection{Linear hyperbolic system}
\label{exls}

\begin{figure}
    \centering
    \includegraphics[width=6.0cm]{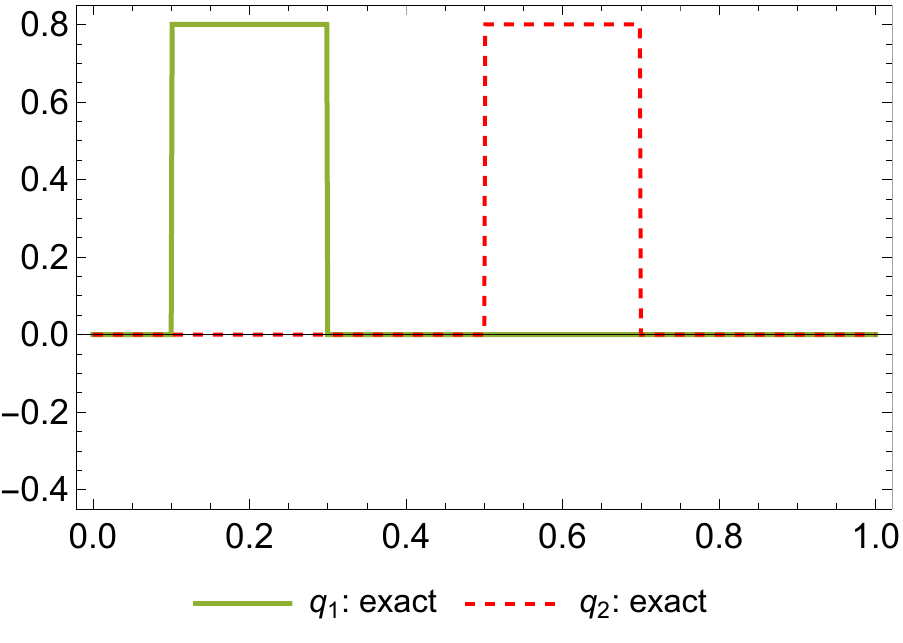}\vspace{.3cm}

    \includegraphics[width=6.0cm]{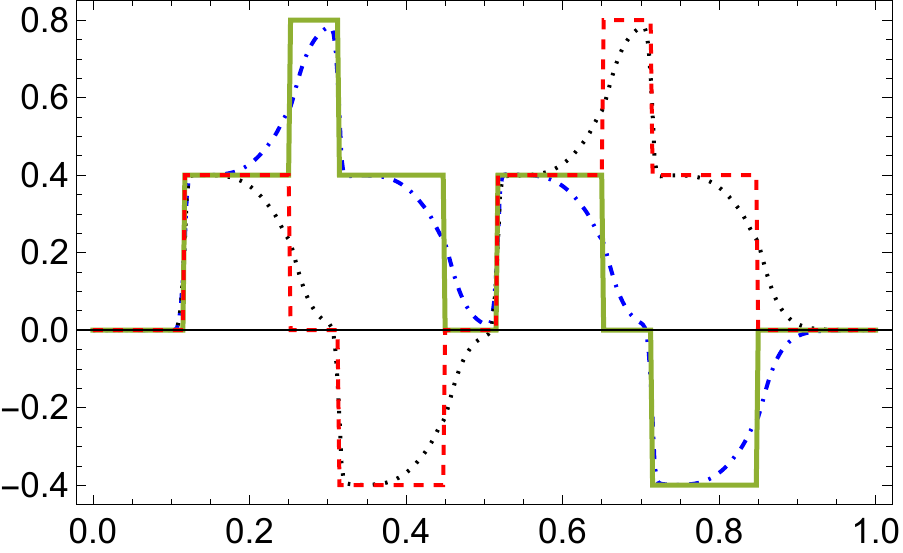}
    \includegraphics[width=6.0cm]{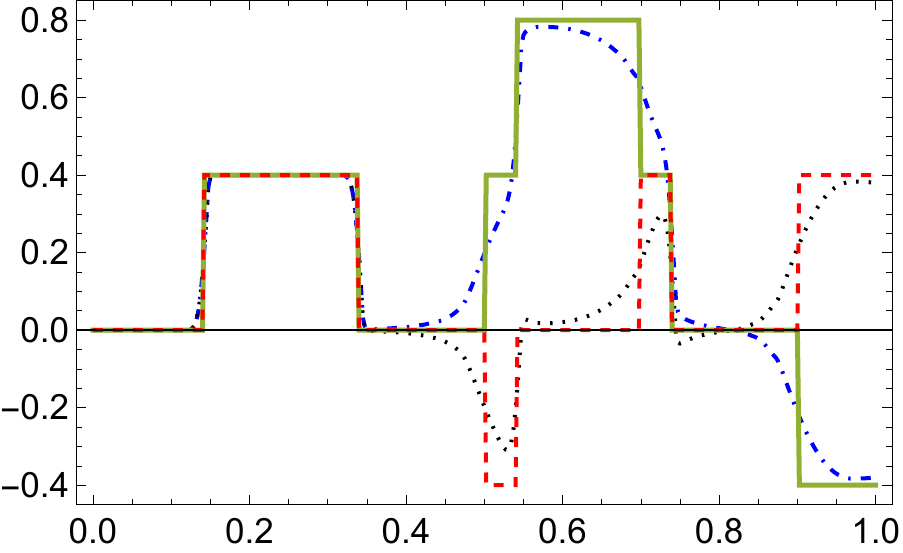}\vspace{.3cm}

    \includegraphics[width=6.0cm]{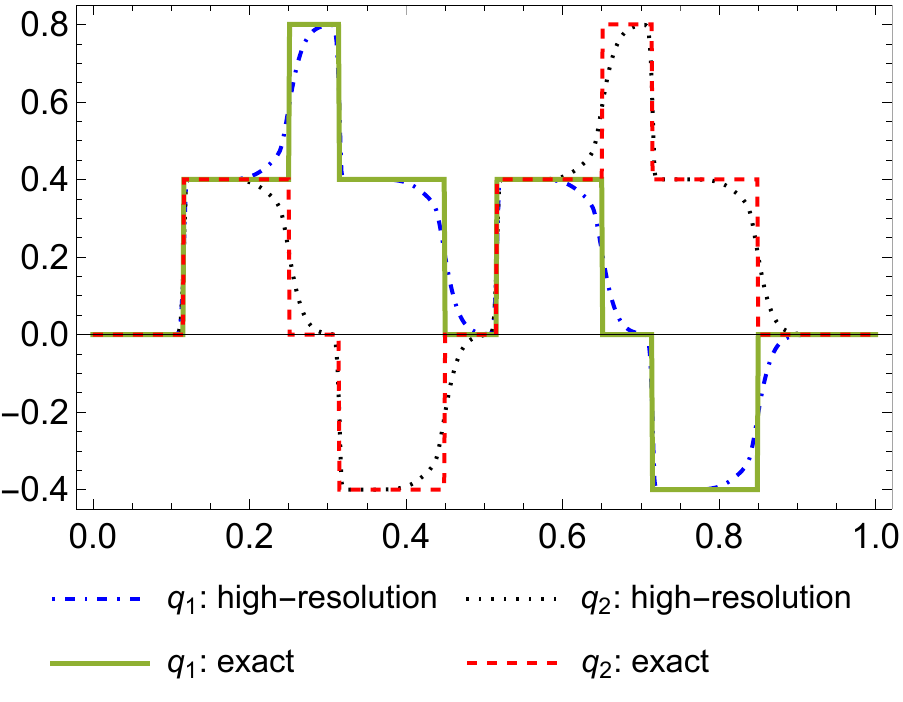}
    \includegraphics[width=6.0cm]{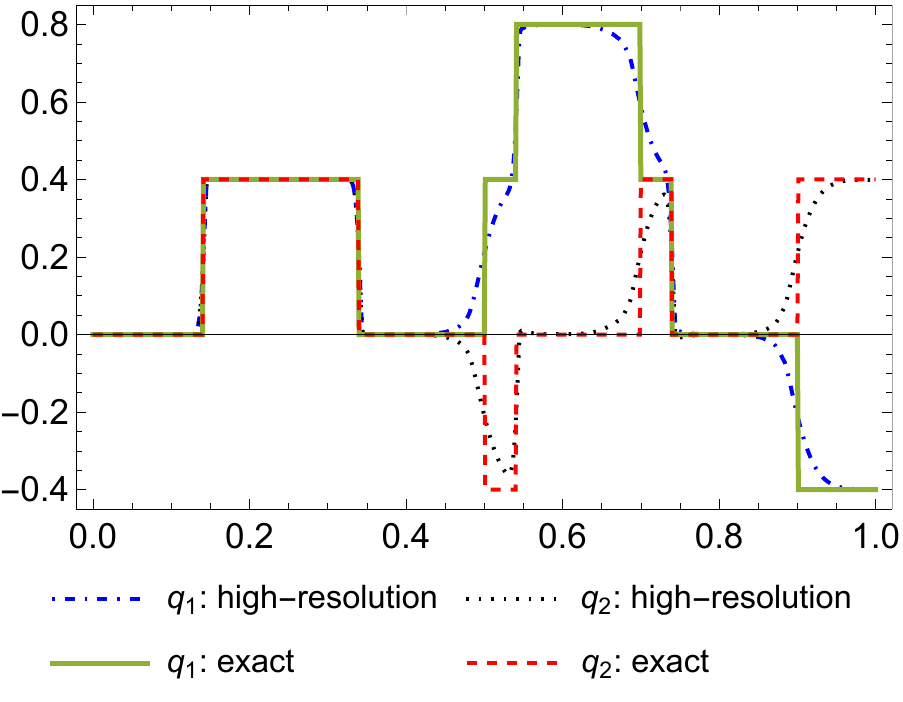}
    
    \caption{The comparison of the exact solutions with the numerical solutions obtained with the high resolution method for the example in Section \ref{exls}. The first row is the initial condition, the second row is for $I=400$ and $t=0.15$ (the first column) and $t=0.4$ (the second column), and the third row is for $I=800$ and the analogous times. The constant Courant number is $10$ using $\tau=10h$.}
    \label{fig:linsys}
\end{figure}

To test the method for systems of conservation laws (\ref{sys}) as described in Section \ref{sec4}, we begin with a simple linear equation having a constant matrix,
\begin{equation}
    \label{exlins}\nonumber
    {\bf f}={\bf f}({\bf q})=\textbf{f}(q_1,q_2) = A \cdot {\bf q} \,, \quad 
    A = \frac{1}{2} \left (\begin{array}{cc}
        1.1 & -0.9 \\\\
        -0.9 & 1.1 
    \end{array}
    \right) \,.
\end{equation}
The matrix $A$ has positive eigenvalues $1$ and $0.1$ that can formally represent the fast and the slow characteristic speed, respectively. The initial functions consist of rectangular profiles
$$
q_1(x,0)=\left\{\begin{array}{lr}
   0.8  & 0.1 < x < 0.3 \\
    0 & \hbox{otherwise}
\end{array}\right. \,, \quad
q_2(x,0)=\left\{\begin{array}{lr}
   0.8  & 0.5 < x < 0.7 \\
    0 & \hbox{otherwise}
\end{array}\right. \,,
$$
see the first row in Figure \ref{fig:linsys}. The problem is considered for $x \in [0,1]$ and the exact solution is defined by
\begin{eqnarray}
\nonumber
q_1(x,t)=\frac{1}{2}\left( q_1(x-0.1t,0)+q_1(x-t,0)+q_2(x-0.1t,0)-q_2(x-t,0)\right) \,, \\[1ex]
\nonumber
q_1(x,t)=\frac{1}{2}\left( q_1(x-0.1t,0)-q_1(x-t,0)+q_2(x-0.1t,0)+q_2(x-t,0)\right)  \,.
\end{eqnarray}

The example is computed with the Courant number $10$ using $\tau=10h$, so only the slowly moving waves can be well resolved with numerical solutions. The results are presented at two different times in Figure \ref{fig:linsys} where one can clearly see that the numerical solutions do not contain any visible oscillations and that the contact discontinuities are well resolved for the slowly moving waves and smeared for the fast moving discontinuities. Clearly, the choice of time step for this example is dictated only by accuracy requirements and not by stability restrictions.

\subsection{Shallow water equations}
\label{exsw1}

\begin{figure}
    \centering
    \includegraphics[width=6.0cm]{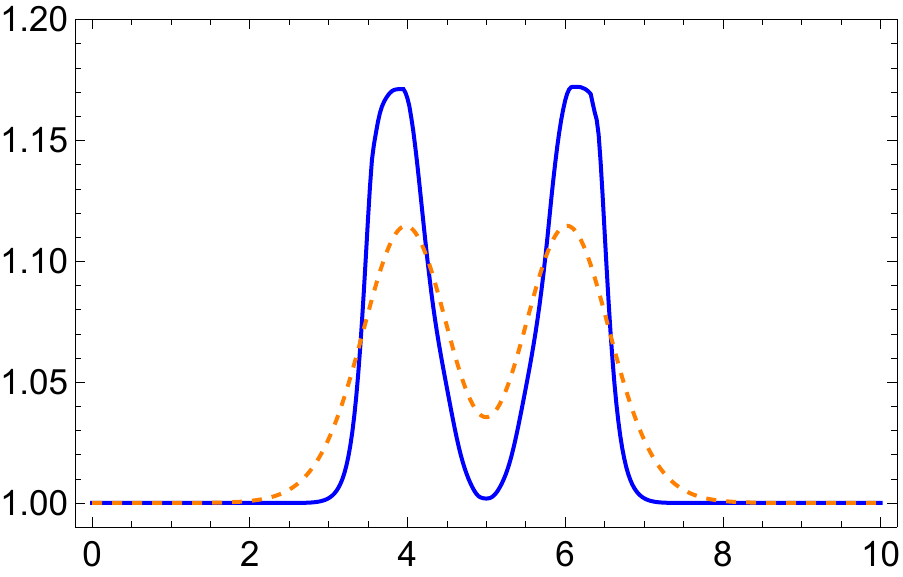}
    \includegraphics[width=6.0cm]{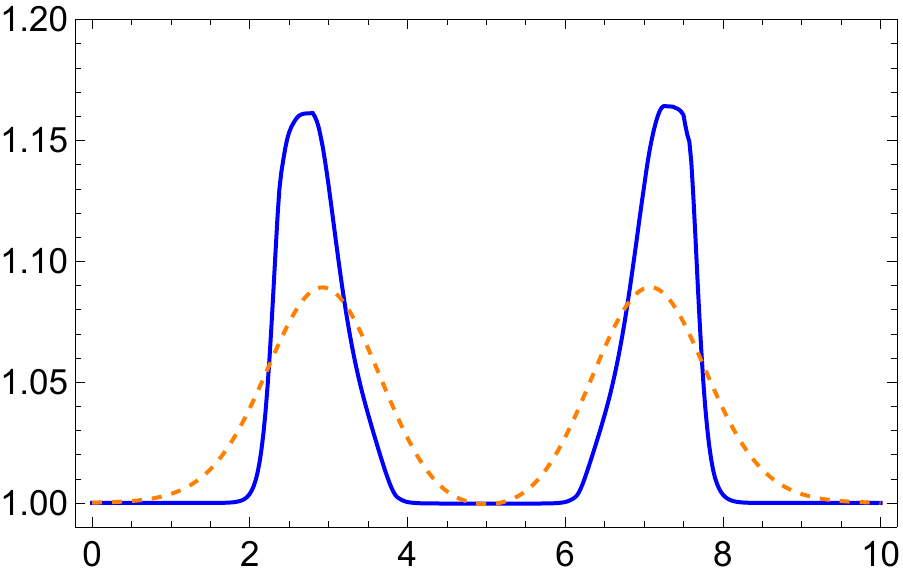}\vspace{.3cm}

    \includegraphics[width=6.0cm]{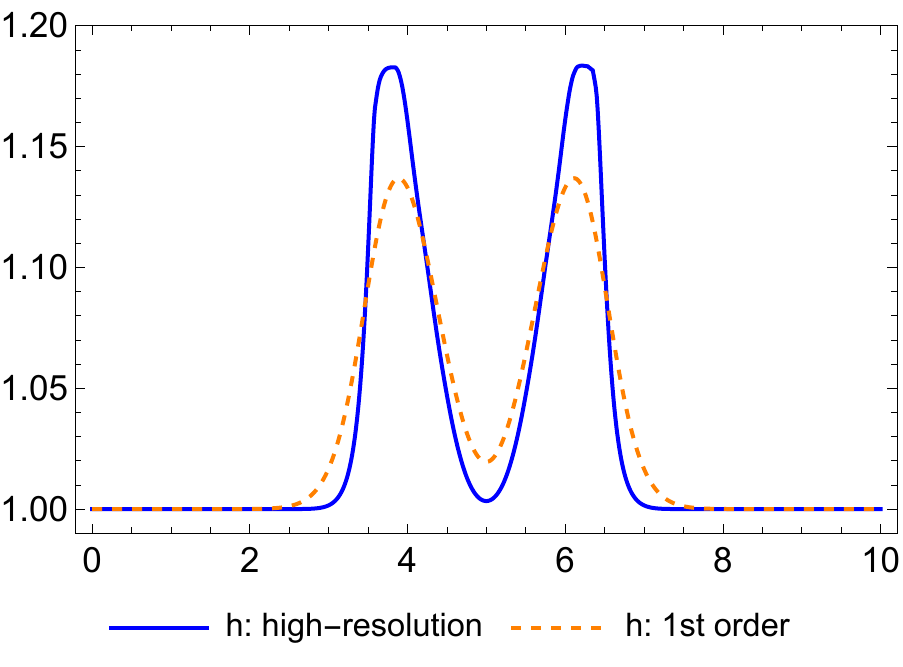}
    \includegraphics[width=6.0cm]{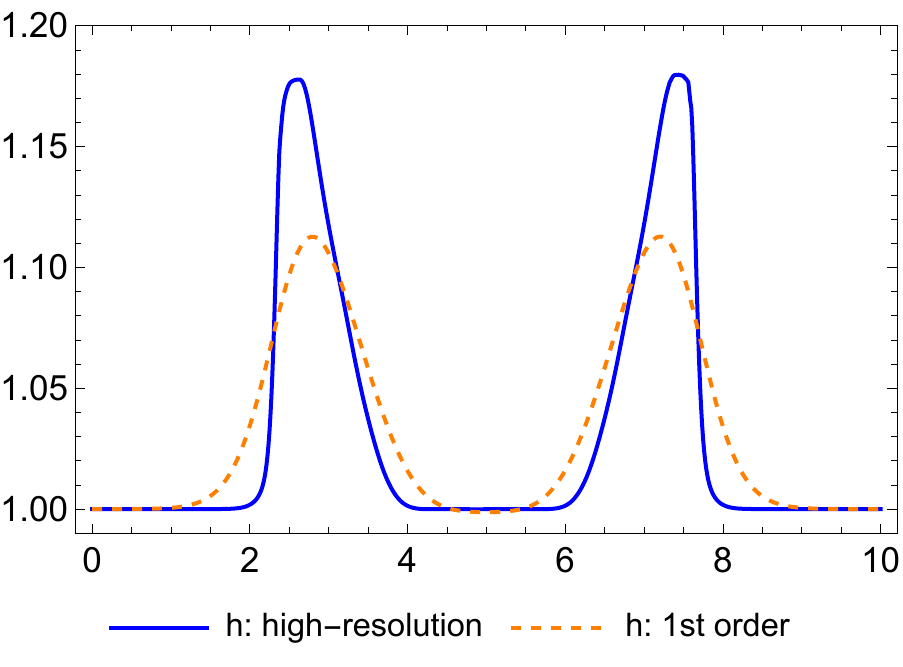}   \vspace{.5cm} 
    
    \includegraphics[width=6.0cm]{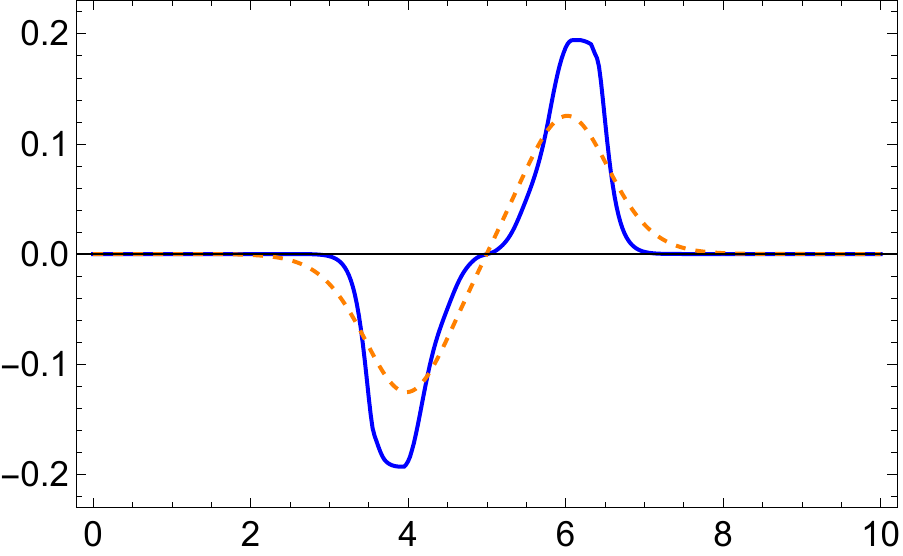}
    \includegraphics[width=6.0cm]{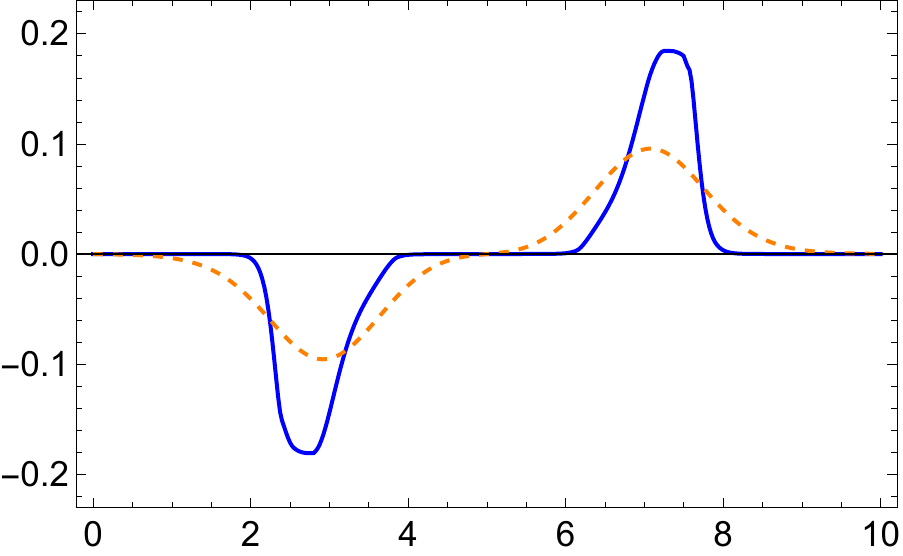}\vspace{.3cm}

    \includegraphics[width=6.0cm]{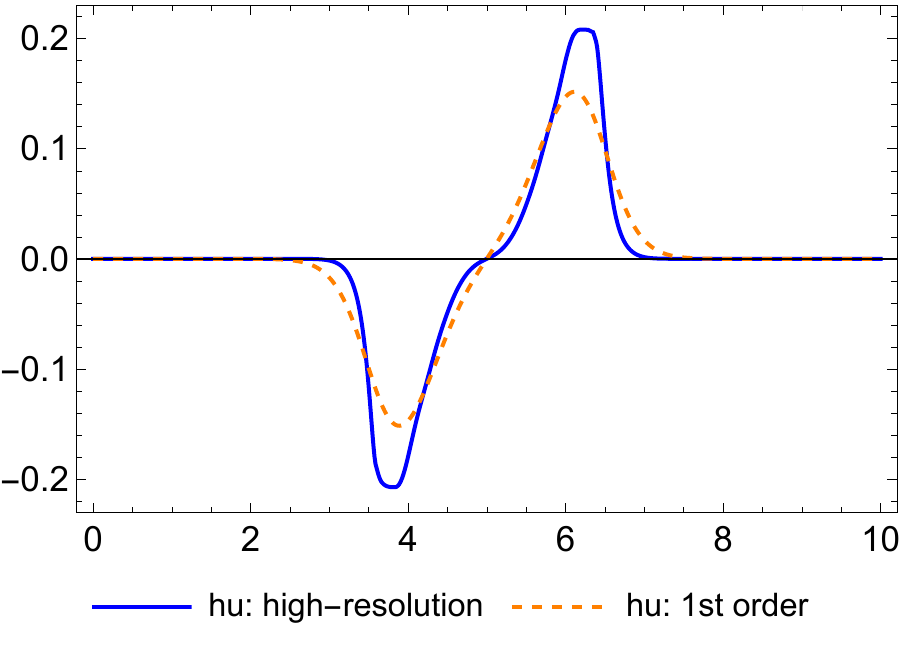}
    \includegraphics[width=6.0cm]{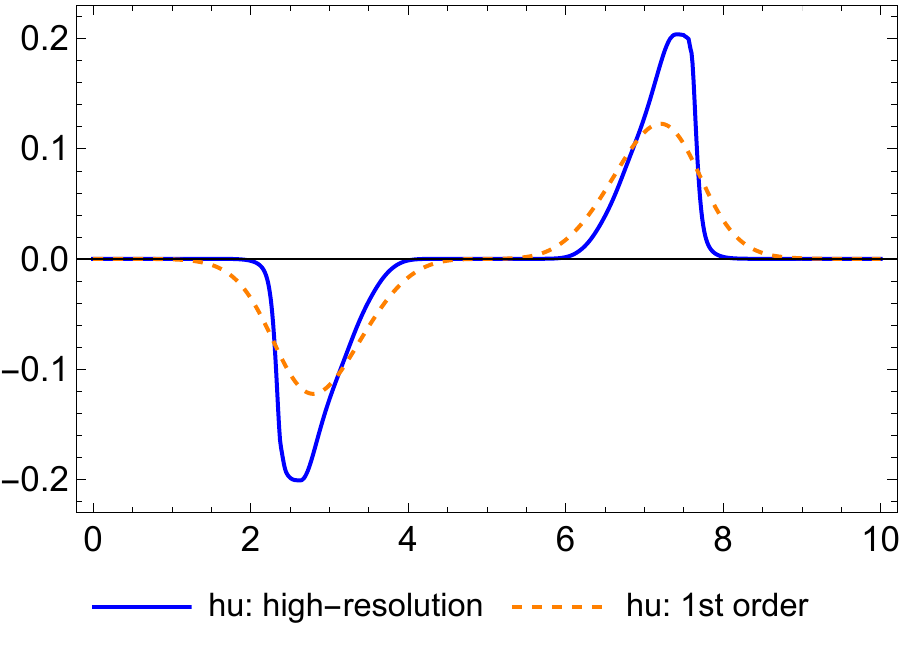}    
    \caption{The comparison of numerical solutions obtained with the first order method and the high resolution method for the example in Section \ref{exsw1}. The first column is for $t=1$, the second one for $t=2$. The first row compares $h$ for $I=400$, the second one $h$ for $I=800$, the third one $h u$ for $I=400$ and the fourth one $h u$ for $I=800$. The maximal Courant number is always $6.21$.}
    \label{fig:shalow1}
\end{figure}

\begin{figure}
    \centering
    \includegraphics[width=6.0cm]{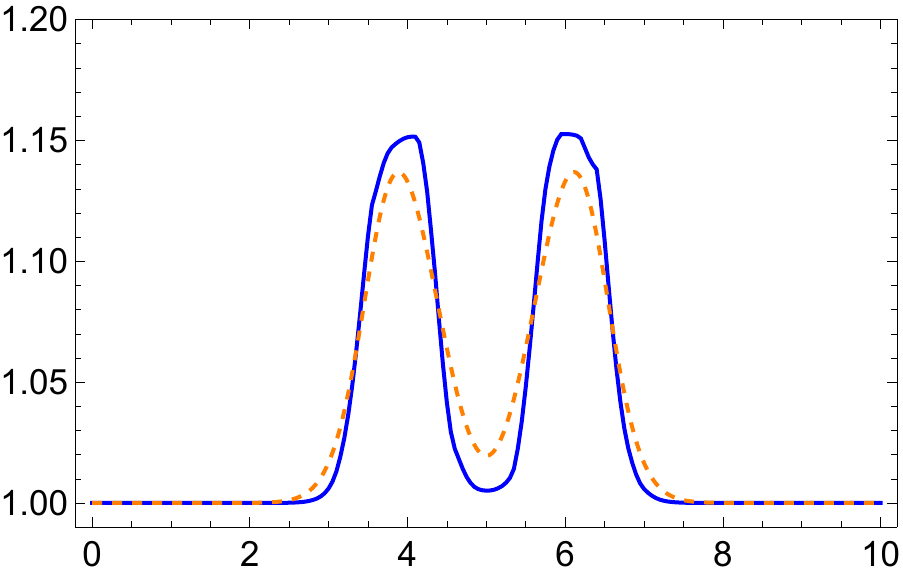}
    \includegraphics[width=6.0cm]{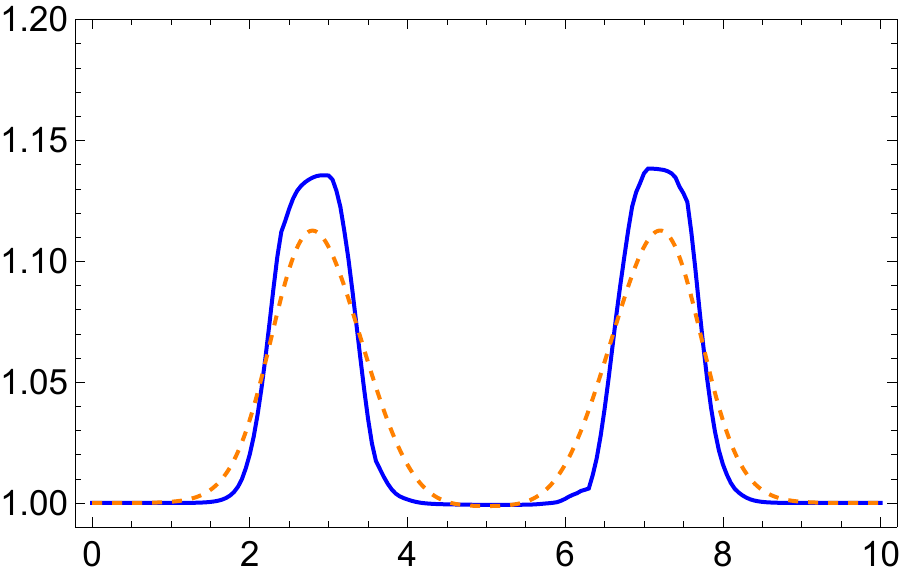}\vspace{.6cm}
    
    \includegraphics[width=6.0cm]{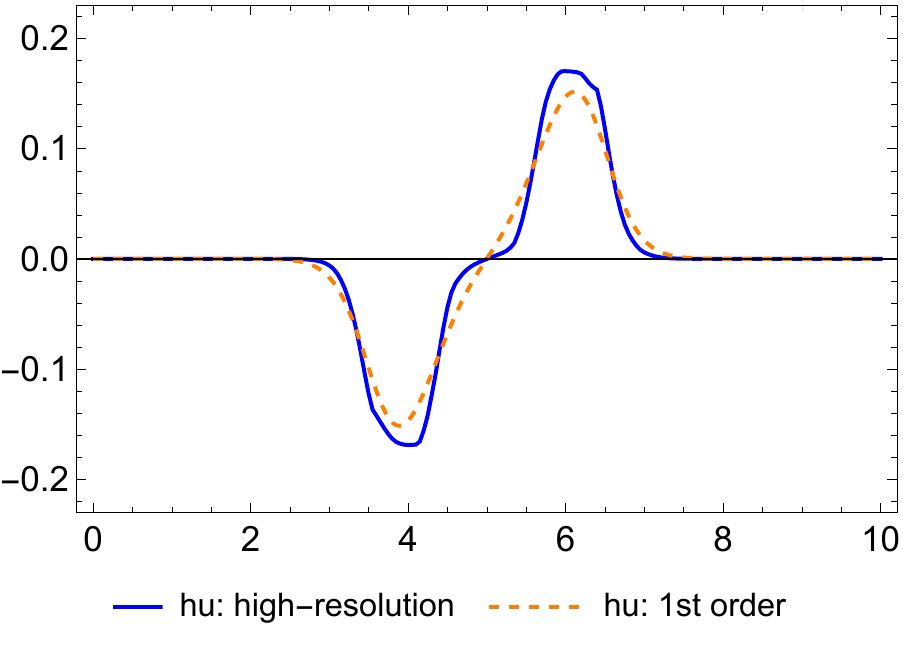}
    \includegraphics[width=6.0cm]{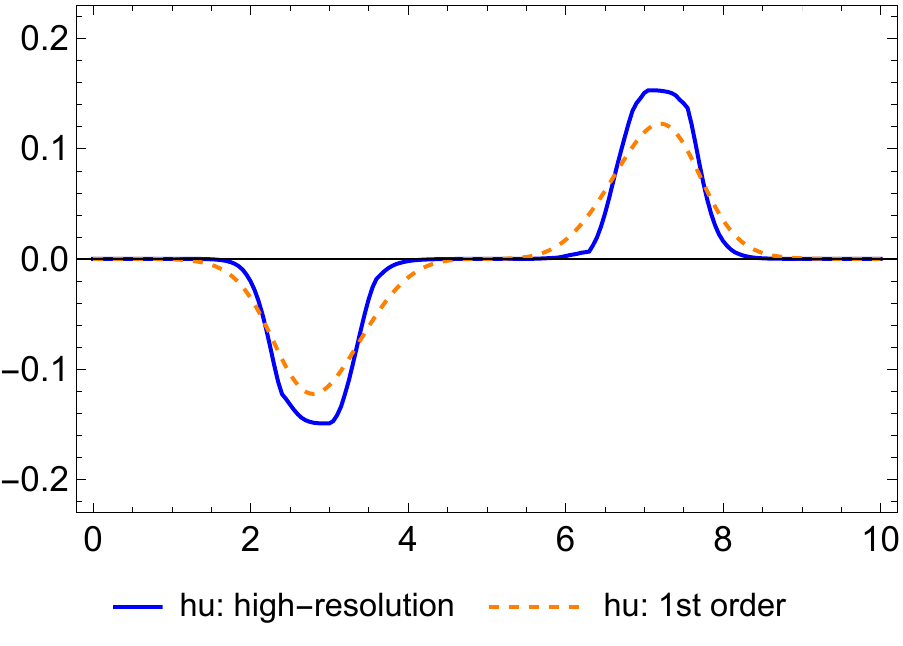}   
    \caption{The comparison of numerical solutions obtained with the first order method and the high resolution method at $t=1$ (left) and $t=2$ (right). The first order method is computed with $I=800$ and the high resolution one with $I=200$. The maximal Courant number is always $6.21$.}
    \label{fig:shalow2}
\end{figure}

Finally, to test the high resolution method in a general case we compute a simple example of the nonlinear hyperbolic system represented by the shallow water equations with the initial condition taken from \cite{leveque_finite_2004},
\begin{eqnarray}
    \label{exsweh}\nonumber
    \partial_t h + \partial_x (h u) = 0 \,, \quad h(x,0) = 1 + 0.4 e^{-5 (x-5)^2} \,,\\[1ex]
    \label{exsweh2}\nonumber
    \partial_t (h u) + \partial_x (h u^2 + 0.5 h^2) = 0 \,, \quad u(x,0) = 0 \,.
\end{eqnarray}
The system is considered for $x \in [0,10]$ and $t \in [0,2]$ and the constant boundary conditions $h=1$ and $u=0$ are used in $x=0$. The system is discretized with conservative variables $(h,h u)$ using the Lax-Friedrichs splitting (\ref{lf}). The eigenvalues and the eigenvectors of ${\bf f}'$ can be analytically expressed \cite{leveque_finite_2004}, so the value of $\alpha$ in (\ref{lf}) is set to $1.3$ by a rough estimate of the maximal eigenvalues for the expected values of $h$ and $u$ to preserve the inequalities in (\ref{split}). Note that, e.g., the choice $\alpha=1.2$ slightly violated the inequalities in (\ref{split}) during computations. 

The comparison of results at $t=1$ and $t=2$ for the first order scheme and the high resolution scheme is given in Figure \ref{fig:shalow1} for two fine grids with $\tau = 5 h$ giving the maximal Courant number around $6.21$. One can see a significantly improved resolution of the shock and the rarefaction waves when comparing the high resolution method with the first order accurate one. The results resemble well those presented in \cite{leveque_finite_2004}.

To make the difference of the resolution even clearer, we compare in Figure \ref{fig:shalow2} the results obtained on a coarse grid with the high resolution method and the results obtained by the first order accurate method on twice uniformly refined grid which still do not have the quality of the high resolution method.

\section{Conclusion}
\label{sec-conc}

We have presented the compact implicit conservative finite difference method for hyperbolic problems in the one-dimensional case. The method shares the advantageous properties of the first order accurate implicit method in \cite{lozano_implicit_2021}. Namely, for the linear advection equation with constant speed, the method is unconditionally stable, and the numerical solutions are obtained explicitly after one (forward or backward) step of the fast sweeping method \cite{frolkovic_semi-implicit_2021}. In the case of nonlinear scalar hyperbolic PDEs, one has to solve for each grid point in each step of the fast sweeping method a single nonlinear algebraic equation with the nonlinearity only due to the nonlinear flux function. All these properties are preserved in the proposed high resolution TVD method which is second order accurate if the solution is smooth. Although the TVD limiters depend on the unknown solution, this nonlinearity can be typically resolved with one predictor and one corrector step. The method is applied successfully for the linear systems of hyperbolic PDE and for the shallow water equations by expressing the second order correction terms in the scheme using the characteristic variables and speeds.

The proposed high resolution compact implicit method can be considered for the problems where fully implicit or explicit-implicit schemes have appeared useful. In addition to accuracy requirements, there are no other restrictions on the choice of time steps for stability reasons. A possible restriction on the time step due to slow or no convergence of the nonlinear algebraic solver is shared with the first order accurate implicit method. We plan to extend the method analogously to \cite{titarev2005weno,puppo_quinpi_2022} with a high order WENO type spatial reconstruction and a Lax-Wendroff type of time discretization \cite{baeza_approximate_2020}. 

\bibliographystyle{plain}
\bibliography{references}
\end{document}